\documentclass[12pt]{article}
\usepackage{fullpage,graphicx,psfrag,amsmath,version}
\usepackage{url,algorithm}
\usepackage[small,bf]{caption}
\usepackage{algorithmicx,algpseudocode,fixltx2e}
\newcommand{\BIT}{\begin{itemize}}
\newcommand{\EIT}{\end{itemize}}
\newcommand{\BEAS}{\begin{eqnarray*}}
\newcommand{\EEAS}{\end{eqnarray*}}
\newcommand{\BEQ}{\begin{equation}}
\newcommand{\EEQ}{\end{equation}}
\newcommand{\eg}{{\it e.g.}}
\newcommand{\ie}{{\it i.e.}}

\newcommand{\ones}{\mathbf 1}
\newcommand{\reals}{{\mbox{\bf R}}}
\newcommand{\integers}{{\mbox{\bf Z}}}
\newcommand{\complex}{{\mbox{\bf C}}}
\newcommand{\symm}{{\mbox{\bf S}}}  

\newcommand{\argmin}{\mbox{\textrm{argmin}}}

\newcommand{\Rank}{\mathop{\bf Rank}}
\newcommand{\Tr}{\mathop{\bf Tr}}
\newcommand{\diag}{\mathop{\bf diag}}

\newcommand{\card}{\mathop{\bf card}}

\algdef{SE}[DOWHILE]{Do}{doWhile}{\algorithmicdo}[1]{\algorithmicwhile\ #1}%
\algnewcommand{\algorithmicgoto}{\textbf{Go to}}%
\algnewcommand{\Goto}[1]{\algorithmicgoto~\ref{#1}}%
\newcommand{\C}{\mathcal C}
\newcommand{\Crelax}{\mathcal C^\mathrm{rlx}}
\newcommand{\Crestrict}{\mathcal C^\mathrm{rstr}}
\newcommand{\Cneighbor}{\mathcal C^\mathrm{ngbr}}

\title{A General System for Heuristic Solution of Convex Problems
over Nonconvex Sets}
\author{Steven Diamond \and Reza Takapoui \and Stephen Boyd}

\excludeversion{incomplete}

\begin{document}
\maketitle
\begin{abstract}
We describe general heuristics to approximately
solve a wide variety of problems with convex objective
and decision variables from a nonconvex set.
The heuristics, which employ convex relaxations, convex restrictions,
local neighbor search methods,
and the alternating direction method of multipliers (ADMM),
require the solution of a modest number of convex problems,
and are meant to apply to general problems, without much tuning.
We describe an implementation of these methods in a package called NCVX,
as an extension of CVXPY, a Python package
for formulating and solving convex optimization problems.
We study several examples of well known nonconvex problems, and show that
our general purpose heuristics are effective in finding approximate
solutions to a wide variety of problems.
\end{abstract}

\section{Introduction}

\subsection{The problem}

We consider the optimization problem
\BEQ\label{e-prob}
\begin{array}{ll}
\mbox{minimize}   & f_0(x, z) \\
\mbox{subject to} & f_i(x, z) \leq 0, \quad i=1,\ldots,m \\
& Ax + Bz = c \\
& z \in \C,
\end{array}
\EEQ
where $x \in \reals^n$ and $z \in \reals^q$
are the decision variables,
$A \in \reals^{p \times n}$, $B \in \reals^{p \times q}$,
$c \in \reals^p$ are problem data, and $\C\subseteq \reals^q$ is compact.
We assume that the objective and inequality constraint functions
$f_0,\ldots,f_m: \reals^n \times \reals^q \to \reals$
are jointly convex in $x$ and $z$.
When the set $\C$ is convex, (\ref{e-prob}) is a convex
optimization problem,
but we are interested here in the case where $\C$ is not convex.
Roughly speaking, the problem~(\ref{e-prob}) is a convex optimization
problem, with some additional nonconvex constraints, $z\in \C$.
We can think of $x$ as the collection of
decision variables that appear only
in convex constraints, and $z$ as the decision variables that
are directly constrained to lie in the (generally) nonconvex
set $\C$.
The set $\C$ is often a Cartesian product, $\C = \C_1 \times \cdots \times
\C_k$, where
$\C_i \subset \reals^{q_i}$ are sets that are simple to describe,
\eg, $\C_i = \{0,1\}$.
We denote the optimal value of the problem~(\ref{e-prob}) as $p^\star$,
with the usual conventions that $p^\star=+\infty$ if the problem
is infeasible, and $p^\star = -\infty$ if the problem is unbounded
below.

\subsection{Special cases}
\paragraph{Mixed-integer convex optimization.}
When $\C = \{0,1\}^q$, the problem~(\ref{e-prob}) is a general
mixed integer convex program, \ie,
a convex optimization problem in which some variables are constrained
to be Boolean.
(Mixed Boolean convex program would be a more accurate
name for such a problem, but `mixed integer' is commonly used.)
It follows that the problem~(\ref{e-prob}) is hard; it includes as
a special case, for example,
the general Boolean satisfaction problem.

\paragraph{Cardinality constrained convex optimization.}
As another broad special case of~(\ref{e-prob}),
consider the case $\C =\{z \in \reals^q \mid \card(z) \leq k,~
\|z\|_\infty \leq M \}$, where $\card(z)$ is the number of nonzero
elements of $z$, and $k$ and $M$ are given.
We call this the general \emph{cardinality-constrained convex problem}.
It arises in many interesting applications, such as regressor selection.

\paragraph{Other special cases.}
As we will see in \S\ref{s-examples},
many (hard) problems can be formulated in the form~(\ref{e-prob}).
More examples include regressor selection, 3-SAT, circle packing,
the traveling salesman problem,
factor analysis modeling, job selection,
the maximum coverage problem, inexact graph isomorphism, and many more.

\subsection{Convex relaxation}
\paragraph{Convex relaxation of a set.}
A compact set $\C$ always has a tractable \emph{convex relaxation}.
By this we mean a (modest-sized) set of convex inequality and linear
equality constraints that hold for every $z\in \C$:
\[
z\in \C ~\Longrightarrow~ h_i(z) \leq 0, \quad i=1,\ldots,s, \qquad Fz=g.
\]
We will assume that these relaxation constraints are included in the convex
constraints of~(\ref{e-prob}).
Adding these relaxation constraints
to the original problem yields an equivalent problem (since the
added constraints are redundant),
but can improve the convergence of any method, global or heuristic.
By tractable, we mean that the number of added constraints is modest,
and in particular, polynomial in $q$.

For example, when $\C = \{0,1\}^q$, we have the inequalities
$0 \leq z_i \leq 1$, $i=1,\ldots, q$.
(These inequalities define the convex hull of $\C$, \ie, all other
convex inequalities that hold for all $z\in \C$ are implied by them.)
When
\[
\C =\{z \in \reals^q \mid \card(z) \leq k,~
\|z\|_\infty \leq M \},
\]
we have the convex inequalities
\[
\|z\|_1 \leq kM, \qquad \|z\|_\infty \leq M.
\]
(These inequalities define the convex hull of $\C$.)
For general compact $\C$ the inequality $\|z\|_\infty \leq M$ will
always be a convex relaxation for some $M$.

\paragraph{Relaxed problem.}
If we remove the nonconvex constraint $z\in\C$,
we get a \emph{convex relaxation} of the original problem:
\BEQ\label{e-relaxation}
\begin{array}{ll}
\mbox{minimize}   & f_0(x,z)\\
\mbox{subject to} & f_i(x, z) \leq 0, \quad i=1,\ldots,m\\
& Ax + Bz = c.
\end{array}
\EEQ
(Recall that convex equalities and inequalities known to hold for
$z\in \C$ have been incorporated in the convex constraints.)
The relaxed problem is convex; its optimal value is a lower bound
on the optimal value $p^\star$ of~(\ref{e-prob}).
A solution $(x^*, z^*)$ to problem~(\ref{e-relaxation}) need not satisfy
$z^* \in \C$, but if it does, the pair $(x^*,z^*)$ is optimal for~(\ref{e-prob}).

\begin{incomplete}
\paragraph{Breaking symmetry.}
In a convex relaxation, we add inequalities that are known to hold
for \emph{any} point $z \in \C$.
When there are symmetries in the problem~(\ref{e-prob}),
we can add additional convex constraints to the problem, with the
property that there is always at least one optimal point that satisfies
the new constraints.  (It follows that solving the problem with the
additional constraints is equivalent to the original problem.)

The symmetry typically arises due to arbitrary labeling of objects,
for example when we are
assigning colors, or machines or jobs from a pool of equivalent ones.
We describe a simple case here to give the idea; we will see symmetry
arising in several examples as well.

Suppose the problem~(\ref{e-prob}) is invariant under permutations of the
variables $w_1, \ldots, w_t \in \reals^\ell$ (which are subvectors of $(x,z)$).
That is, permuting the $w_i$ does not change the objective,
or change their feasibility.   Let $h\in \reals^\ell$ be arbitrary.
We can without loss of generality (by permuting $w_i$)
assume that $h^T w_i$ are sorted, \ie,
\[
h^T w_1 \leq \cdots \leq h^T w_t.
\]
Thus we can impose this constraint on the problem without loss of
generality.
As a variation, we can add the (convex) penalty function
\[
\tau(x,z) = \sum_{i=1}^{t-1} \max\{h^T(w_i - w_{i+1}),0\},
\]
which penalizes non-sorted choices of $w_i$,
to the objective.
It is easy to see that this yields an equivalent problem:
given any feasible point $(x,z)$ for~(\ref{e-prob}),
we can always permute $w_i$ so that $\tau(x,z) = 0$.

Unlike the basic convex relaxations, which pertain to the set $\C$
and are independent of the convex parts of the problem,
these symmetry-breaking constraints or penalty depend on the
entire problem, and not just the set $\C$.
\end{incomplete}

\subsection{Projections and approximate projections}
Our methods will make use of tractable projection, or
tractable approximate projection, onto the set $\C$.  The usual Euclidean
projection onto $\C$ will be denoted $\Pi$.
(It need not be unique when $\C$ is not convex.)
By approximate projection,
we mean any function $\hat \Pi: \reals^q \to \C$ that satisfies
$\hat \Pi (z) = z $ for $z \in \C$.
For example, when $\C = \{0,1\}^q$, exact projection is given by
rounding the entries to $\{0,1\}$.

As a less trivial example, consider the cardinality-constrained problem.
The projection of $z$ onto $\C$ is given by
\[
\left(\Pi \left(z\right)\right)_i = \left\{ \begin{array}{ll}
M & z_i > M, ~i \in \mathcal I\\
-M & z_i < -M, ~i \in \mathcal I\\
z_i & |z_i| \leq M, ~i \in \mathcal I\\
0 & i \not\in \mathcal I,
\end{array}\right.
\]
where $\mathcal I \subseteq \{1,\ldots, q\}$ is a set of indices of
$k$ largest values of $|z_i|$.
We will describe many projections, and some approximate projections,
in \S\ref{s-set}.

\subsection{Residual and merit functions}
For any $(x,z)$ with $z \in \C$, we define the \emph{constraint
residual} as
\[
r(x,z) = \sum_{i=1}^m (f_i(x,z))_+ + \|Ax+Bz-c\|_1,
\]
where $(u)_+=\max\{u,0\}$ denotes the positive part;
$(x,z)$ is feasible if and only if $r(x,z)=0$.
Note that $r(x,z)$ is a convex function of $(x,z)$.
We define the \emph{merit function} of a pair $(x,z)$ as
\[
\eta (x,z)= f_0(x,z)+\lambda r(x,z),
\]
where $\lambda>0$ is a parameter.
The merit function is also a convex function of $(x,z)$.

When $\C$ is convex and the problem is feasible,
minimizing $\eta(x,z)$ for large enough $\lambda$
yields a solution of the original problem~(\ref{e-prob})
(that is, the residual is a so-called exct penalty function);
when the problem is not feasible, it tends to find
approximate solutions that satisfy many of the constraints
\cite{han1979exact,di1989exact,fletcher1973exact}.

We will use the merit function to judge candidate approximate solutions
$(x,z)$ with $z \in \C$; that is, we take a pair with lower
merit function value to be a better approximate solution than one
with higher merit function value.
For some problems (for example, unconstrained problems)
it is easy to find feasible points, so all candidate points
will be feasible.
The merit function then reduces to the objective value.
At the other extreme, for feasibility problems the objective is
zero, and goal is to find a feasible point.  In this case the
merit function reduces to $\lambda r(x,z)$, \ie, a positive
multiple of the residual function.

\subsection{Solution methods}

In this section we describe various methods for solving the
problem~(\ref{e-prob}), either exactly (globally) or approximately.

\paragraph{Global methods}
Depending on the set $\C$, the problem~(\ref{e-prob}) can be solved
globally by a variety of algorithms, including (or mixing)
branch-and-bound
\cite{lawler1966branch,narendra1977branch,brucker1994branch},
branch-and-cut
\cite{padberg1991branch,tawarmalani2005polyhedral,stubbs1999branch},
semidefinite hierarchies \cite{sherali1990hierarchy}, or even
direct enumeration when $\C$ is a finite set.
In each iteration of these methods, a convex optimization problem derived
from~(\ref{e-prob}) is solved, with $\C$ removed, and (possibly) additional
variables and convex constraints added.
These global methods are generally thought to have high worst-case
complexities and indeed can be very slow in practice, even for modest
size problem instances.

\paragraph{Local solution methods and heuristics}
A local method for~(\ref{e-prob}) solves a modest number of convex
problems, in an attempt to find a good approximate solution,
\ie, a pair $(x,z)$ with $z\in \C$ and a low value of the merit
function $\eta(x,z)$.
For a feasibility problem, we might hope to find a solution;
and if not, find one with a small constraint residual.
For a general problem, we can hope to find a feasible point with
low objective value, ideally near the lower bound on $p^\star$ from
the relaxed problem.
If we cannot find any feasible points, we can settle for a pair $(x,z)$
with $z \in \C$ and low merit function value.
All of these methods are heuristics, in the sense that they cannot
in general be guaranteed to find an optimal, or even good, or even
feasible, point in only a modest number of iterations.

There are of course many heuristics for the general
problem~(\ref{e-prob}) and for many of its special cases.
For example, any global optimization method can be stopped
after some modest number of iterations; we take the best point found
(in terms of the merit function) as our approximate solution.
We will discuss some local search methods, including neighbor search and
polishing, in \S\ref{local}.

\paragraph{Existing solvers}
There are numerous open source and commercial solvers that can handle
problems with nonconvex constraints. We only mention a few of them here.
Gurobi \cite{gurobi}, CPLEX \cite{cplex2009v12}, MOSEK \cite{mosek}
provide global methods for mixed integer linear programs,
mixed integer quadratic programs, and mixed integer
second order cone programs.
BARON \cite{ts:05}, Couenne \cite{cplex2009v12}, and
SCIP \cite{achterberg2009scip} use global methods for
nonlinear programs
and mixed integer nonlinear programs.
Bonmin \cite{bonami2008algorithmic} and Knitro \cite{byrd2006knitro}
provide global methods for mixed integer convex programs and
heuristic methods for mixed integer nonliner programs.
IPOPT \cite{waechter2009introduction} and
NLopt \cite{johnson2014nlopt} use heuristic methods for nonlinear
programs.

\subsection{Our approach}
The purpose of this paper is to describe a general system
for heuristic solution of~(\ref{e-prob}), based on solving
a \emph{modest} number of convex problems derived from~(\ref{e-prob}).
By heuristic, we mean that the algorithm need not find an optimal point,
or indeed, even a feasible point, even when one exists.
We would hope that for many feasible problem instances from some
application, the algorithm does find
a feasible point, and one with objective not too far from
the optimal value.
The disadvantage of a heuristic over a global method is clear and simple:
it need not find an optimal point.
The advantage of a heuristic is that it can be (and often is)
dramatically faster to carry out than a global method.
Moreover there are many applications where a heuristic method
for~(\ref{e-prob}) is sufficient.  This might be the case when
the objective and constraints are already approximations of what
we really want, so the added effort of solving it globally is not
worth it.

\paragraph{ADMM.}
One of the heuristic methods described in this paper
is based on the alternating directions method of multipliers
(ADMM), an operator splitting algorithm originally devised to solve
convex optimization problems \cite{boyd2011distributed}.
We call this heuristic nonconvex alternating directions method of
multipliers  (NC-ADMM).
The idea of using ADMM as a heuristic to solve nonconvex problems
was mentioned in \cite[Ch. 9]{boyd2011distributed}, and has been explored by
Yedidia and others \cite{Derbinsky2013AnImp} as a message
passing algorithm. Consensus ADMM has been used
for general quadratically constrained quadratic programming
in \cite{huang2016consensus}. In \cite{xu2012alternating}, ADMM
has been applied to non-negative matrix factorization with missing
values. ADMM also has been used for real and complex polynomial
optimization models in \cite{jiang2014alternating},for constrained
tensor factorization in \cite{liavas2014parallel}, and for optimal
power flow in \cite{erseghe2014distributed}.
ADMM is a generalization of the method of multipliers
\cite{hestenes1969multiplier, bertsekas2014constrained},
and there is a long history of using the method of multipliers
to (attempt to) solve nonconvex problems
\cite{chartrand2012nonconvex, chartrand2013nonconvex,hong2014distributed,
hong2014convergence,peng2015proximal,wang2014convergence,LP:15}.
Several related methods, such as the Douglas-Rachford method
\cite{eckstein1992douglas}
or Spingarn's method of partial inverses \cite{spingarn1985applications},
could just as well have been used.

\paragraph{Our contribution.}
The paper has the following structure.
In \S\ref{s-algorithm} we discuss local search methods and describe how
they can be used as solution improvement methods. This will enable us to
study simple but sophisticated methods such as relax-round-polish, and
iterative neighbor search.
In \S\ref{s-set} we catalog a variety of nonconvex
sets for which Euclidean projection or approximated projection
is easily evaluated and, when applicable, we discuss relaxations, restrictions,
and the set of neighbors for a given point.
In \S\ref{s-implementation} we discuss an implementation
of our general system for heuristic solution NCVX,
as an extension of CVXPY \cite{cvxpy},
a Python package for
formulating and solving convex optimization problems.
The object-oriented features of CVXPY make the extension
particularly simple to implement.
Finally, in \S\ref{s-examples} we demonstrate the performance
of our methods on several example problems.

\section{Local improvement methods}
\label{local}
In this section we describe some simple general local search methods.
These methods take a point $z\in\C$ and by performing a local search on $z$ they
find a candidate pair $(\hat x, \hat z)$, with $\hat z\in\C$ and a
lower merit function. We will see that for many applications
these methods with a good initialization can be used to obtain an
approximate solution. We will also see how we can use these methods
to improve solution candidates from other heuristics, hence we refer to
these methods as \emph{solution improvement}.

\subsection{Polishing}
\paragraph{Convex restriction.}
We can have a tractable \emph{convex restriction} of $\C$,
that includes a given point in $\C$.
This means that for each point $\tilde z\in\C$,
we have a set of convex equalities and inequalities on $z$, that hold for
$\tilde z$, and imply $z \in \C$.
We denote the set of points that satisfy the restrictions as
$\Crestrict (\tilde z)$, and call this set the \emph{restriction} of
$\C$ at $\tilde z$.  The restriction set $\Crestrict(\tilde z)$ is convex,
and satisfies $\tilde z \in \Crestrict (\tilde z) \subseteq \C$.
The trivial restriction is given by $\Crestrict(\tilde z)=\{\tilde z\}$.

When $\C$ is discrete, for example $\C= \{0,1\}^q$, the trivial restriction
is the only restriction. In other cases we can have interesting nontrivial
restrictions, as we will see below.
For example, with
$\C=\{z\in \reals^q \mid \card(z)\leq k,~ \|z\|_\infty \leq M\}$,
we can take as restriction
$\Crestrict(\tilde z)$, the set of vectors $z$ with
the same sparsity pattern as $\tilde z$, and $\|z \|_\infty \leq M$.

\paragraph{Polishing.}
Given any point $\tilde z\in\C$, we can replace the constraint $z\in\C$
with $z\in \Crestrict(\tilde z)$ to get the convex problem
\BEQ\label{e-polish}
\begin{array}{ll}
\mbox{minimize}   &\eta(x,z)\\
\mbox{subject to}  &z\in \Crestrict(\tilde z),
\end{array}
\EEQ
with variables $x,z$.
(When the restriction $\Crestrict(\tilde z)$ is the trivial one, \ie,
a singleton, this is equivalent to fixing $z=\tilde z$ and
minimizing over $x$.)
We call this problem the \emph{convex restriction} of~(\ref{e-prob})
at the point $\tilde z$.
The restricted problem is convex,
and its optimal value is an upper bound on $p^\star$.

As a simple example of polishing consider the mixed
integer convex problem. The only restriction is the trivial one, so
the polishing problem for a given Boolean vector $\tilde z$
simply fixes the values of the Boolean variables, and solves the
convex problem over the remaining variables, \ie, $x$.
For the cardinality-constrained convex problem, polishing
fixes the sparsity pattern of $z$ and solves the resulting convex
problem over $z$ and $x$.

For problems with nontrivial restrictions, we can solve the
polishing problem repeatedly until convergence. In other
words we can use the output of the polishing problem as
an initial point for another polishing problem and keep iterating
until convergence or until a maximum number of iterations is reached.
This technique is called \emph{iterated polishing} and
described in algorithm \ref{iter-polish}.
\begin{algorithm}[H]
\caption{Iterated polishing}
\begin{algorithmic}[1]
\State Input: $\tilde z$
\Do
\State $z^\mathrm{old} \gets \tilde z$.
\State Find $(\tilde x, \tilde z)$
by solving the polishing problem with restriction
$z \in \Crestrict(z^\mathrm{old})$.
\doWhile{$\tilde z \neq z^\mathrm{old}$}\\
\Return $(\hat x,\hat z)$.
\end{algorithmic}
\label{iter-polish}
\end{algorithm}

If there exists a point $\tilde x$ such that $(\tilde x, \tilde z)$ is feasible, the
restricted problem is feasible too.
The restricted problem need not be feasible in general, but if it is,
with solution $(\hat x, \hat z)$, then the pair
$(\hat x,\hat z)$ is feasible for the original problem~(\ref{e-prob})
and satisfies $f(\hat x,\hat z) \leq f(\tilde x,\tilde z)$
for any $\tilde x$ for which $(\tilde x,\tilde z)$ is feasible.
So polishing can take a point $\tilde z\in\C$ (or a pair
$(\tilde x,\tilde z)$) and produce another pair $(\hat x, \hat z)$ with
a possibly better objective value.

\subsection{Relax-round-polish}
With the simple tools described so far
(\ie, relaxation, polishing, and
projection) we can create several heuristics for approximately
solving the problem~(\ref{e-prob}).
A basic version solves the relaxation, projects the relaxed
value of $z$ onto $\C$, and then polishes the result.

\begin{minipage}[c]{1 \textwidth}
\begin{algorithm}[H]
\caption{Relax-round-polish heuristic}
\begin{algorithmic}[1]
\State Solve the convex relaxation~(\ref{e-relaxation}) to obtain
$(x^\mathrm{rlx}, z^\mathrm{rlx})$.
\State Find $z^\mathrm{rnd}=\Pi(z^\mathrm{rlx})$.
\State Find $(\hat x, \hat z)$
by solving the polishing problem with restriction
$z \in \Crestrict(z^\mathrm{rnd})$.
\end{algorithmic}
\label{alg_summary}
\end{algorithm}
\end{minipage}
\bigskip

Note that in the first step we also obtain a lower bound on the
optimal value $p^\star$; in the polishing step we obtain
an upper bound, and a feasible pair $(\hat{x},\hat{z})$ that achieves the upper
bound (provided that polishing is successful).
The best outcome is for these bounds to be equal, which means that
we have found a (global) solution of~(\ref{e-prob}) (for this problem
instance).
But relax-round-polish can fail; for example, it can fail to find a
feasible point even though one exists.

Many variations on relax-round-polish are possible.
We can introduce randomization by replacing the round step with
\[
z^\mathrm{rnd}=\Pi(z^\mathrm{rlx}+w),
\]
where $w$ is a random vector.
We can repeat this heuristic with $K$ different random instances of $w$.
For each of $K$ samples of $w$, we polish, giving us a set of $K$ candidate
approximate solutions.
We then take as our final approximate solution the best among
these $K$ candidates,
\ie,
the one with least merit function.

\subsection{Neighbor search} \label{s-neighbor}
\paragraph{Neighbors.}
We describe the concept of \emph{neighbors} for a point $z\in\C$ when
$\C$ is discrete. The set of neighbors of a point $z\in\C$,
denoted $\Cneighbor(z)$, is the set of points
with distance one from $z$ in a natural (integer valued) distance,
which depends on the set $\C$.
For example for the set of Boolean vectors in $\reals^n$ we use
\emph{Hamming distance},
the number of entries in which two Boolean vectors differ.
Hence neighbors of a Boolean
vector $z$ are the set of vectors that differ from $z$ in one component.
The distance between two permutation matrices is defined as the minimum number of
swaps of adjacent rows and columns necessary to transform the first
matrix into the second.
With this distance,
neighbors of a permutation matrix $Z$ are the
set of permutation matrices generated by swapping any two adjacent rows
or columns in $Z$.

For Cartesian products of discrete sets we use the sum of
distances.  In this case, for
$z = (z_1,z_2,\ldots,z_k) \in \C = \C_1 \times \C_2 \times \ldots \times \C_k$, neighbors
of $z$ are points of the form
$(z_1,\ldots,z_{i-1},\tilde z_i,z_{i+1},\ldots, z_k)$ where
$\tilde z_i$ is a neighbor of $z_i$ in $\C_i$.

\paragraph{Basic neighbor search.}
We introduced polishing as a tool that can find a pair $(\hat x,\hat z)$ given an
input $\tilde z\in\C$ by solving a sequence of convex problems. In basic
neighbor search we solve the polishing problem for $\tilde z$ and all neighbors
of $\tilde z$ and return the pair $(x^*, z^*)$ with the smallest merit function value.
In practice, we can sample from $\Cneighbor(\tilde z)$ instead of iterating over
all points in $\Cneighbor(\tilde z)$ if $\left|\Cneighbor(\tilde z)\right|$ is large.

\begin{algorithm}[H]
\caption{Basic neighbor search}
\begin{algorithmic}[1]
\State Input: $\tilde z$
\State Initialize $(x_\mathrm{best}, z_\mathrm{best})=\emptyset$, $\eta_\mathrm{best}=\infty$.
\For{$\hat z \in \{ \tilde z \} \cup \Cneighbor(\tilde z)$}
\State Find $(x^*, z^*)$,
by solving the polishing problem~(\ref{e-polish}),
with constraint $z\in \Crestrict(\hat z)$.
\If{$\eta(x^*, z^*) < \eta_\mathrm{best}$}
\State $(x_\mathrm{best}, z_\mathrm{best})= (x^*, z^*)$, $\eta_\mathrm{best} =\eta(x^*, z^*)$.
\EndIf
\EndFor\\
\Return $(x_\mathrm{best},z_\mathrm{best})$.
\end{algorithmic}
\label{iter-hill-climb}
\end{algorithm}

\paragraph{Iterated neighbor search.}
We can carry out the described neighbor search iteratively as follows.
We maintain a current value of $z$, corresponding to the
best pair $(x,z)$ found so far.  We then consider a neighbor of $z$ and
polish.  If the new point $(x,z)$ is better than the current best one,
we reset our best and continue; otherwise we examine another neighbor.
This is done until a maximum number of iterations is reached, or
all neighbors of the current best $z$ produce (under polishing) no better
pairs. This procedure is sometimes called \emph{hill climbing},
since it resembles an attempt to find the top of a mountain by repeatedly
taking steps towards an ascent direction.

\begin{algorithm}[H]
\caption{Iterative neighbor search}
\begin{algorithmic}[1]
\State Input: $\tilde z$
\State Find $(x_\mathrm{best}, z_\mathrm{best})$ by solving the polishing problem~(\ref{e-polish})
\State$\eta_\mathrm{best}\gets\eta(x_\mathrm{best}, z_\mathrm{best})$.
\For{$\hat z \in\Cneighbor(z_\mathrm{best})$}\label{marker}
\State Find $(x^*, z^*)$,
by solving the polishing problem~(\ref{e-polish}),
with constraint $z\in \Crestrict(\hat z)$.
\If{$\eta(x^*, z^*) < \eta_\mathrm{best}$}
\State $(x_\mathrm{best}, z_\mathrm{best})= (x^*, z^*)$, $\eta_\mathrm{best} =\eta(x^*, z^*)$.
\State \Goto{marker}
\EndIf
\EndFor\\
\Return $(x_\mathrm{best},z_\mathrm{best})$.
\end{algorithmic}
\label{alg:iter-neighbor}
\end{algorithm}
Notice that when no neighbors are available for $\tilde z\in\C$,
this algorithm reduces to simple polishing.

\section{NC-ADMM}\label{s-algorithm}
We already can use the simple tools described in the previous section
as heuristics to find approximate solutions to problem~(\ref{e-prob}).
In this section, we describe the alternating direction method of multipliers (ADMM)
as a mechanism to generate candidate points $\tilde z$ to carry out local search methods
such as iterated neighbor search. We call this method nonconvex ADMM, or NC-ADMM.
\subsection{ADMM}
Define $\phi:\reals^q \to \reals \cup \{ -\infty, +\infty\}$ such that
$\phi(z)$ is the best objective value of problem~(\ref{e-prob})
after fixing $z$. In other words,
\[
\phi(z) = \inf_x \left\{ f_0(x,z) \mid
f_i(x, z) \leq 0,~i=1,\ldots,m,~ Ax + Bz = c \right\}.
\]
Notice that $\phi(z)$ can be $+\infty$ or $-\infty$ in case the problem is not feasible
for this particular value of $z$, or
problem~(\ref{e-relaxation}) is unbounded below after fixing $z$.
The function $\phi$ is convex, since it is the
partial minimization of a convex function over a
convex set \cite[\S 3.4.4]{boyd2004convex}.
It is defined over all points $z\in\reals^q$, but we are
interested in finding its minimum value over
the nonconvex set $\C$. In other words,
problem~(\ref{e-prob}) can be formulated as
\BEQ\label{e-prob-ref}
\begin{array}{ll}
\mbox{minimize}   & \phi(z) \\
\mbox{subject to} & z\in\C.
\end{array}
\EEQ

As discussed in \cite[Chapter\ 9]{boyd2011distributed}, ADMM can be used
as a heuristic to solve nonconvex constrained problems. ADMM has the form
\begin{equation}
\label{e-update}
  \begin{split}
w^{k+1} &:=  \argmin_z \left(\phi(z)+ (\rho/2)\|z - z^k + u^k\|_2^2\right)\\
z^{k+1} &:=  \Pi\left(w^{k+1} - z^{k} + u^k\right)\\
u^{k+1} &:=  u^k + w^{k+1} - z^{k+1},
\end{split}
\end{equation}
where $\rho>0$ is an algorithm parameter,
$k$ is the iteration counter,
and $\Pi$ denotes Euclidean projection onto $\C$
(which need not be unique when $\C$ is not convex).

The initial values $u^0$ and $z^0$ are additional algorithm parameters.
We always set $u^0 = 0$ and draw $z^0$
randomly from a normal distribution $\mathcal N(0, \sigma^2 I)$,
where $\sigma > 0$ is an algorithm parameter.

\subsection{Algorithm subroutines}
\paragraph{Convex proximal step}
Carrying out the first step of the
algorithm, \ie, evaluating the proximal operator of $\phi$,
involves solving the convex optimization problem
\BEQ\label{e-prox-step}
\begin{array}{ll}
\mbox{minimize}   & f_0(x,z) + (\rho/2) \|z-z^k+u^k\|_2^2 \\
\mbox{subject to} & f_i(x, z) \leq 0, \quad i=1,\ldots,m,\\
& Ax + Bz = c,
\end{array}
\EEQ
over the variables $x \in \reals^n$ and $z \in \reals^q$.
This is the original problem~(\ref{e-prob}), with the nonconvex
constraint $z\in \C$ removed, and an additional
convex quadratic term involving $z$ added to the objective.
We let $(x^{k+1}, w^{k+1})$ denote a solution of~(\ref{e-prox-step}).
If the problem~(\ref{e-prox-step}) is infeasible, then so is
the original problem~(\ref{e-prob}); should this happen, we can terminate
the algorithm with the certain conclusion that~(\ref{e-prob})
is infeasible.

\paragraph{Projection}
The (nonconvex) projection step consists of finding the closest point in $\C$ to
$w^{k+1} - z^{k} + u^k$.
If more than one point has the smallest distance,
we can choose one of the minimizers arbitrarily.

\paragraph{Dual update}
The iterate $u^k \in \reals^q$
can be interpreted as a scaled dual variable, or
as the running sum of the error values $w^{k+1}-z^k$.

\subsection{Discussion}
\paragraph{Convergence.}
When $\C$ is convex (and a solution of~(\ref{e-prob}) exists),
this algorithm is guaranteed to converge to a solution,
in the sense that $f_0(x^{k+1},w^{k+1})$ converges to the optimal
value of the problem~(\ref{e-prob}),
and $w^{k+1}-z^{k+1}\to 0$, \ie, $w^{k+1}\to \C$.
See~\cite[\S 3]{boyd2011distributed} and the references therein for a
more technical description and details.
But in the general case, when $\C$ is not convex,
the algorithm is not guaranteed to converge, and even when it does,
it need not be to a global, or even local, minimum.
Some recent progress has been made on understanding convergence
in the nonconvex case \cite{LP:15}.

\paragraph{Parameters.}
Another difference with the convex case is that the convergence and the
quality of solution depends on $\rho$, whereas for convex problems
this algorithm is guaranteed to converge to the optimal value regardless
of the choice of $\rho$.
In other words, in the convex case the choice of parameter $\rho$ only
affects the speed of the convergence, while in the nonconvex case the
choice of $\rho$ can have a
critical role in the quality of approximate solution,
as well as the speed at which this solution is found.

The optimal parameter selection for ADMM is still an active research area.
In \cite{ghadimi2015optimal} the optimal parameter selection
for quadratic problems is discussed.
In a more generalized setting, Giselsson discusses the optimal parameter selection
for ADMM for strongly convex functions \cite{giselsson2014diagonal,
giselsson2014monotonicity,
giselsson2014preconditioning}.
The dependency of global and local convergence
properties of ADMM on parameter choice has been studied in \cite{hong2012linear,
boley2013local}.

\paragraph{Initialization.}
In the convex case the choice of initial point $z^0$ affects the number
of iterations to find a solution, but not the quality of the solution.
Unsurprisingly, the nonconvex case differs in that the choice of $z^0$
has a major effect on the the quality of the approximate solution.
As with the choice of $\rho$,
the initialization in the nonconvex case is currently an active
area of research; see, \eg, \cite{huang2016consensus,LP:15,takapoui2015simple}.
Getting the best possible results on a particular problem requires
a careful and problem specific choice of initialization.
We draw initial points randomly from $\mathcal N(0,\sigma^2 I)$
because we want a method that generalizes easily across many different problems.

\subsection{Solution improvement}
Now we describe two techniques to obtain better solutions after carrying
out ADMM. The first technique relies on iterated neighbor search and
the second one is using multiple restarts with random initial points in order to
increase the chance of obtaining a better solution.

\paragraph{Iterated neighbor search}
After each iteration, we can carry out iterated polishing (as described in \S \ref{s-neighbor})
with $\Crestrict(z^{k+1})$ to obtain $(\hat x^{k+1},\hat z^{k+1})$. We will return the pair
with the smallest merit function as the output
of the algorithm.

\paragraph{Multiple restarts}
As we mentioned, we choose the initial value $z^0$ from a normal
distribution $\mathcal N(0,\sigma^2 I)$. We can run the algorithm multiple
times from different initial points to increase the chance of a feasible point
with a smaller objective value.

\subsection{Overall algorithm}
The following is a summary of the algorithm with solution improvement.

\begin{algorithm}[H]
\caption{NC-ADMM heuristic} 
\begin{algorithmic}[1]
\State Initialize $u^0=0$, $(x_{\textrm{best}},z_{\textrm{best}}) = \emptyset$, $\eta_\mathrm{best}=\infty$.
\For{algorithm repeats $1, 2, \ldots, M$}
\State Initialize $z^0\sim\mathcal N(0,\sigma^2I)$.
\For{$k=1, 2, \dots, N$}
\State $(x^{k+1}, w^{k+1}) \gets \argmin_z \left(\phi(z)+ (\rho/2)\|z - z^k + u^k\|_2^2\right)$.
\State $z^{k+1} \gets  \Pi\left(w^{k+1} - z^{k} + u^k \right)$.
\State Use algorithm (\ref{alg:iter-neighbor}) on $z^{k+1}$ to get the improved iterate $(\hat x, \hat z)$.
\If {$\eta(\hat x, \hat z)<\eta_{\textrm{best}}$}
\State $(x_{\textrm{best}},z_{\textrm{best}}) \gets (\hat x, \hat z)$, $\eta_{\textrm{best}} = \eta(\hat x, \hat z)$.
\EndIf
\State $u^{k+1} \gets  u^k + w^{k+1} - z^{k+1}$.
\EndFor
\EndFor \\
\Return $x_{\textrm{best}},z_{\textrm{best}}$.
\end{algorithmic}
\label{alg_summary}
\end{algorithm}
\bigskip

\section{Projections onto nonconvex sets}\label{s-set}
In this section we catalog various nonconvex sets
with their implied convex constraints which will be included in the convex
constraints of
problem~(\ref{e-prob}). We also provide a Euclidean
projection (or approximate projection)
$\Pi$ for these sets.
Also, when applicable, we introduce a nontrivial restriction and set
of neighbors.

\subsection{Subsets of $\reals$}
\paragraph{Booleans}
For $\C = \{0,1\}$, a convex relaxation (in fact, the convex hull of $\C$) is $[0,1]$.
Also, a projection is simple rounding:
$\Pi(z) = 0$ for $z \leq 1/2$, and $\Pi(z) =1$ for $z > 1/2$.
($z=1/2$ can be mapped to either point.)
Moreover, $\Cneighbor(0)=\{1\}$ and $\Cneighbor(1)=\{0\}$.

\paragraph{Finite sets}
If $\C$ has $M$ elements, the convex hull of $\C$ is the interval from the
smallest to the largest element. We can project onto $\C$ with
no more than $\log_2 M$ comparisons.
For each $z\in\C$ the set of neighbors of $\C$ are the
immediate points to the right and left of $z$ (if they exist).

\paragraph{Bounded integers}
Let $\C = \integers \cap [-M,M]$, where $M > 0$.
The convex hull is the interval from the smallest to the largest element
integer in $[-M,M]$,
\ie, $[-\lfloor M \rfloor, \lfloor M \rfloor]$.
The projection onto $\C$ is simple: if $z>\lfloor M \rfloor$ ($z< - \lfloor M \rfloor$)
then $\Pi(z) = \lfloor M \rfloor$
($\Pi(z) = -\lfloor M \rfloor$).
Otherwise, the projection of $z$ can be found by simple rounding.
For each $z \in \C$ the set of neighbors of $\C$
is $\{z-1,z+1\}\cap[-M,M]$

\subsection{Subsets of $\reals^n$}

\paragraph{Boolean vectors with fixed cardinality}\label{sec:bool_card}
Let $\C = \{ z\in \{0,1\}^n \mid \card(z) = k\}$.
Any $z \in \C$ satisfies $0\leq z\leq 1$ and $\ones^T z=k$.
We can project $z \in \reals^n$ onto $\C$
by setting the $k$ entries of $z$ with largest value to
one and the remaining entries to zero.
For any point $z\in\C$, the set of neighbors of $z$ is all points
generated by swapping an adjacent $1$ and $0$ in $z$.

\paragraph{Vectors with bounded cardinality}
Let $\C=\{x \in [-M,M]^n \mid \card(x) \leq k\}$, where $M > 0$ and
$k \in \integers_+$.
(Vectors $z\in \C$ are called $k$-sparse.)
Any point $z\in\C$ satisfies
$-M\leq z\leq M$ and $-Mk\leq \ones^T z\leq Mk$.
The projection $\Pi(z)$ is found as follows
\[
\left(\Pi \left(z\right)\right)_i = \left\{ \begin{array}{ll}
M & z_i > M, ~i \in \mathcal I\\
-M & z_i < -M, ~i \in \mathcal I\\
z_i & |z_i| \leq M, ~i \in \mathcal I\\
0 & i \not\in \mathcal I,
\end{array}\right.
\]
where $\mathcal I \subseteq \{1,\ldots, n\}$ is a set of indices of
$k$ largest values of $|z_i|$.

A restriction of $\C$ at $z\in\C$ is the set of
all points in $[-M,M]^n$ that have the same sparsity pattern as $z$.
For any point $z\in\C$, the set of neighbors of $z$ are all points
$x \in \C$ whose sparsity pattern $\tilde{x} \in \{0,1\}^n$ is a neighbor
of $z$'s sparsity pattern $\tilde{z} \in \{0,1\}^n$.
In other words, $\tilde{x}$ can be obtained by swapping an adjacent
$1$ and $0$ in $\tilde{z}$.

\paragraph{Quadratic sets}
Let $\symm_{+}^n$ and $\symm_{++}^n$ denote the set of $n\times n$
symmetric positive semidefinite and symmetric positive definite matrices,
respectively.
Consider the set
\[
\C = \{ z\in \reals^n \mid \alpha \leq z^T A z + 2b^Tz \leq \beta \},
\]
where $A \in \symm_{++}^n$,
$b \in \reals^n$, and $\beta \geq \alpha \geq -b^TA^{-1}b$.
We assume $\alpha \geq -b^TA^{-1}b$ because $z^T A z + 2b^Tz  \geq -b^TA^{-1}b$
for all $z \in \reals^n$.
Any point $z\in\C$ satisfies the convex inequality $z^TAz + 2b^Tz \leq \beta$.

We can find the projection onto $\C$ as follows.
If $z^T A z + 2b^Tz > \beta$, it suffices to solve
\BEQ
\begin{array}{ll}
\mbox{minimize}   & \|x - z\|_2^2 \\
\mbox{subject to} & x^T A x + 2b^Tx \leq \beta,
\end{array}
\EEQ
and if $z^T A z + 2b^Tz < \alpha$, it suffices to solve
\BEQ
\begin{array}{ll}
\mbox{minimize}   & \|x - z\|_2^2 \\
\mbox{subject to} & x^T A x + 2b^Tx \geq \alpha.
\end{array}
\EEQ
(If $\alpha \leq z^T A z + 2b^Tz \leq \beta$, clearly $\Pi(z)=z$.)
The first problem is a convex quadratically constrained quadratic program
and the second problem can be solved by
solving a simple semidefinite program as described in \cite[Appendix B]{boyd2004convex}.
Furthermore, there is a more
efficient way to find the projection by finding the roots of a single-variable polynomial of
degree $2p+1$, where $p$ is the number of distinct eigenvalues of $A$
\cite{huang2016consensus,hmam2010quadratic}.
Note that the projection can be easily found even if $A$ is not positive definite;
we assume $A\in\symm_{++}^n$ only to make $\C$ compact and have a useful convex relaxation.

A restriction of $\C$ at $z \in \C$ is the set
\[
\Crestrict(z) = \{x \in \reals^n \mid
\frac{x^TAz + b^T(x + z) + b^TA^{-1}b}{\sqrt{z^TAz + 2b^Tz + b^TA^{-1}b}} \geq
\sqrt{\alpha + b^TA^{-1}b},
\quad x^TAx + 2b^Tx \leq \beta \}.
\]
Recall that $z^TAz + 2b^Tz + b^TA^{-1}b \geq 0$ for all $z \in \reals^n$
and we assume $\alpha \geq -b^TA^{-1}b$,
so $\Crestrict(z)$ is always well defined.

\paragraph{Annulus and sphere}
Consider the set
\[
\mathcal C = \{z \in \reals^n \mid r \leq \|z\|_2 \leq R\},
\]
where $R \geq r$.

Any point $z\in\C$ satisfies
$\|z\|_2\leq R$.
We can project $z \in \reals^n \setminus \{0\}$
onto $\C$ by the following scaling
\[
\Pi(z) = \begin{cases}
rz/\|z\|_2 & \mbox{if $\|z\|_2 < r$}\\
z & \mbox{if $z\in\C$}\\
R z/\|z\|_2 & \mbox{if $\|z\|_2 > R$},
\end{cases}
\]
If $z = 0$, any point with Euclidean norm $r$ is a valid projection.

A restriction of $\C$ at $z\in\C$ is the set
\[
\Crestrict(z) = \{x \in \reals^n \mid x^Tz\geq r\|z\|_2,~ \|x\|_2 \leq R\}.
\]
Notice that if $r=R$, then $\C$ is a sphere and the restriction will be a singleton.

\paragraph{Box complement and cube surface}
Consider the set
\[
\mathcal C = \{z \in \reals^n \mid a \leq \|z\|_\infty \leq b\}.
\]
Any point $z\in\C$ satisfies $\|z\|_\infty \leq b$. For any point $z$ we can find
the projection $\Pi(z)$ by projecting $z$ component-wise onto $[a,b]$.

Given $z\in\C$ we can obtain a restriction by finding $k=\argmin_i \max\{|z_i|,a\}$ and if $z_k \geq 0$
then
\[
\Crestrict(z) = \{x \mid x_k\geq a,~ \|x\|_\infty \leq b\}.
\]
If $z_k<0$, then
\[
\Crestrict(z) = \{x \mid x_k \leq -a,~ \|x\|_\infty \leq b\}.
\]
Notice that if $a = b$, then $\C$ is a cube surface.


\subsection{Subsets of $\reals^{m\times n}$}
Remember that the projection of a point $X\in\reals^{m \times n}$ on a set
$\C\subset\reals^{m \times n}$ is a point $Z\in\C$ such that
the Frobenius norm $\|X-Z\|_\mathrm F$ is minimized.
As always, if there is more
than one point $Z$ that minimizes
$\|X-Z\|_\mathrm F$, we accept any of them.

\paragraph{Matrices with bounded singular values and orthogonal matrices}
Consider the set of $m \times n$ matrices whose singular values lie between
$1$ and $\alpha$
\[
\C=\{Z\in \reals^{m \times n} \mid I\preceq Z^TZ \preceq \alpha^2 I\},
\]
where $\alpha\geq 1$, and $A \preceq B$ means $B-A\in\symm_+^n$ .
Any point $Z\in\C$ satisfies $\|Z\|_2\leq \alpha$.

If $Z=U\Sigma V^T$ is the singular value decomposition of
$Z$ with singular values $(\sigma_z)_{\min\{m,n\}}\leq\cdots\leq (\sigma_z)_1$ and $X\in\C$
with singular values $(\sigma_x)_{\min\{m,n\}}\leq\cdots\leq (\sigma_x)_1$,
according to the von Neumann trace inequality \cite{von1937some} we will have
\[
\Tr (Z^TX) \leq \sum_{i=1}^{\min\{m,n\}} (\sigma_z)_i(\sigma_x)_i.
\]
Hence \[
\|Z-X\|_F^2 \geq \sum_{i=1}^{\min\{m,n\}} \left((\sigma_z)_i - (\sigma_x)_i\right)^2,
\]
with equality when $X=U\diag(\sigma_x)V^T$.
This inequality implies that $\Pi(Z)=U\tilde\Sigma V^T$, where $\tilde\Sigma$ is a
diagonal matrix and $\tilde\Sigma_{ii}$ is the projection of
$\Sigma_{ii}$ on interval $[1,\alpha]$.
When $Z = 0$, the projection $\Pi(Z)$ is any matrix.

Given $Z=U\Sigma V^T\in\C$, we can have the following restriction \cite{boyd2014mimo}
\[
\Crestrict(Z) = \{X \in \reals^{m \times n} \mid \|X\|_2\leq \alpha,~ V^TX^TU+U^TXV\succeq 2I\}.
\]
(Notice that $X\in\Crestrict(Z)$ satisfies $X^TX\succeq I + (X-UV^T)^T(X-UV^T)\succeq I$.)

There are several noteworthy special cases.
When $\alpha=1$ and $m=n$ we have the set of orthogonal matrices.
In this case, the restriction will be a singleton.
When $n=1$, the set $\C$ is equivalent to the annulus
$\{z \in \reals^m \mid 1 \leq \|z\|_2 \leq \alpha\}$.

\paragraph{Matrices with bounded rank}
Let $\C = \{ Z \in \reals^{m \times n} \mid \Rank(Z) \leq k,~ \|Z\|_2\leq M \}$.
Any point $Z\in\C$ satisfies $\|Z\|_2\leq M$ and
$\|Z\|_*\leq Mk$, where $\|\cdot\|_*$ denotes
the trace norm.
If $Z = U\Sigma V^T$ is the singular value decomposition of $Z$,
we will have $\Pi(Z) = U \tilde \Sigma V^T$, where $\tilde \Sigma$ is
a diagonal matrix with $\tilde \Sigma_{ii} = \min\{\Sigma_{ii},M\}$ for
$i=1,\ldots k$, and $\tilde \Sigma_{ii} =0$ otherwise.

Given a point $Z\in\C$, we can write the singular value decomposition of
$Z$ as $Z=U\Sigma V^T$ with $U\in\reals^{m\times k}$,
$\Sigma\in\reals^{r\times r}$ and $V\in\reals^{n\times k}$.
A restriction of $\C$ at $Z$ is
\[
\Crestrict(Z) = \{U\tilde\Sigma V^T \mid \tilde\Sigma\in\reals^{r\times r}\}.
\]

\paragraph{Assignment and permutation matrices}
The set of \emph{assignment matrices} are Boolean matrices
with exactly one non-zero
element in each column and at most one non-zero element in each row. (They
represent an assignment of the columns to the rows.) In other words, the
set of assignment matrices on $\{0,1\}^{m \times n}$, where $m\geq n$, satisfy
\[
\begin{array}{cc}
\sum_{j=1}^n Z_{ij} \leq 1, &\quad i=1,\ldots,m \\
\sum_{i=1}^m Z_{ij} = 1, &\quad j=1,\ldots,n.
\end{array}
\]
These two sets of inequalities, along with $0\leq Z_{ij}\leq 1$ are the implied
convex inequalities.
When $m=n$, this set becomes the set of permutation matrices, which we show
by $\mathcal P_n$.

Projecting $Z\in\reals^{m\times n}$ (with $m\geq n$) onto the set of
assignment matrices involves choosing an entry from each column of $Z$
such that no two chosen entries are from the same row and the sum of
chosen entries is maximized. Assuming that the entries of $Z$ are the weights
of edges in a bipartite graph, the projection onto the set of assignment matrices will
be equivalent to finding a maximum-weight matching in a bipartite graph.
The Hungarian method \cite{kuhn2005hungarian}
is a well-know polynomial time algorithm
to find the maximum weight matching, and hence also the projection onto assignment
matrices.

The neighbors of an assignment or permutation matrix $Z \in\reals^{m\times n}$
are the matrices generated by swapping two adjacent rows or columns of $Z$.

\paragraph{Hamiltonian cycles}
A \emph{Hamiltonian cycle} is a cycle in a graph that visits every node
exactly once. Every Hamiltonian cycle in a complete graph can be represented
by its adjacency matrix, for example
\[
\left[\begin{array}{cccc}
0 & 0 & 1 & 1       \\
0 & 0 & 1 & 1       \\
1 & 1 & 0 & 0       \\
1 & 1 & 0 & 0
\end{array}\right]
\]
represents a Hamiltonian cycle that visits nodes $(3,2,4,1)$ sequentially.
Let $\mathcal H_n$ be the set of $n\times n$ matrices that represent a
Hamiltonian cycle.

Every point $Z\in\mathcal H_n$ satisfies $0\leq Z_{ij} \leq 1$ for
$i,j = 1,\ldots,n$, and $Z=Z^T$, $(1/2)Z\ones = \ones$, and
\[
2\mathbf I - Z + 4\frac{\ones \ones^T}{n} \geq 2 (1-\cos \frac{2\pi}{n}) \mathbf I,
\]
where $\mathbf I$ denotes the identity matrix. In order to see why the last
inequality holds, it's enough to notice that $2\mathbf I - Z$ is the Laplacian
of the cycle represented by $Z$
\cite{merris1994laplacian, anderson1985eigenvalues}.
It can be shown that the smallest eigenvalue
of $2\mathbf I - Z$ is zero (which corresponds to the eigenvector $\ones$),
and the second smallest eigenvalue of $2\mathbf I - Z$ is
$2 (1-\cos \frac{2\pi}{n})$. Hence all eigenvalues of
$2\mathbf I - Z + 4\frac{\ones \ones^T}{n}$ must be no smaller than
$2 (1-\cos \frac{2\pi}{n})$.

We are not aware of a polynomial time algorithm to find the projection of
a given real $n\times n$ matrix onto $\mathcal H_n$.
We can find an \emph{approximate projection} of $Z$ by the following greedy
algorithm: construct a graph with $n$ vertices where the edge between
$i$ and $j$ is weighted by $z_{ij}$. Start with the edge with largest weight and at
each step, among all the edges that don't create a cycle, choose the edge with
the largest weight (except for the last step where a cycle is created).

For a matrix $Z\in\mathcal H_n$, the set of neighbors of $Z$ are matrices
obtained after swapping two adjacent nodes, \ie, matrices in form $P_{(i,j)}ZP_{(i,j)}^T$
where $Z_{ij}=1$ and
$P_{(i,j)}$ is a permutation matrix that transposes connected
nodes $i$ and $j$ and keeps other nodes unchanged.

\subsection{Combinations of sets}

\paragraph{Cartesian product.}
Let $\C = \C_1 \times \cdots \times \C_k \subset \reals^n$,
where $\C_1,\ldots,\C_k$ are compact sets with known projections (or approximate projections).
A convex relaxation of $\C$ is the Cartesian product
$\Crelax_1 \times \cdots \times \Crelax_k$,
where $\Crelax_i$ is the set described by the convex relaxation of $\C_i$.
The projection of $z \in \reals^n$ onto $\C$ is
$\left(\Pi_1(z_1),\ldots,\Pi_k(z_k)\right)$, where $\Pi_i$ denotes the projection
onto $\C_i$ for $i=1,\ldots,k$.

A restriction of $\C$ at a point
$z = (z_1,z_2,\ldots,z_k) \in \C$
is the Cartesian product
$\Crestrict(z) = \Crestrict_1(z_1) \times \cdots \times \Crestrict_k(z_k)$.
The neighbors of $z$ are all points
$(z_1,\ldots,z_{i-1},\tilde z_i,z_{i+1},\ldots, z_k)$ where
$\tilde z_i$ is a neighbor of $z_i$ in $\C_i$.

\paragraph{Union.}
Let $\C = \cup_{i=1}^k\C_i$,
where $\C_1,\ldots,\C_k$ are compact sets with known projections
(or approximate projections).
A convex relaxation of $\C$ is the constraints
\[
\begin{array}{cc}
x_i \in \Crelax_i,& \quad i = 1,\ldots,k \\
s_i \in [0,1],& \quad i = 1,\ldots,k \\
z = \sum_{i=1}^k x_i& \\
\sum_{i=1}^k s_i = 1& \\
\|x_i\|_\infty \leq M_i s_i,& \quad i = 1,\ldots,k,
\end{array}
\]
where $\Crelax_i$ is the set described by the convex relaxation of $\C_i$
and $M_i > 0$ is the minimum value such that $\|z_i\|_\infty \leq M_i$
holds for all $z_i \in \C_i$.

We can project $z \in \reals^n$ onto $\C$ by projecting onto each set
separately and keeping the projection closest to $z$:
\[
\Pi(z) = \argmin_{x\in\{\Pi_1(z),\cdots,\Pi_k(z)\}}\|z - x\|_2.
\]
Here $\Pi_i$ denotes the projection onto $\C_i$.

A restriction of $\C$ at a point $z$ is $\Crestrict_i(z)$ for any
$\C_i$ containing $z$.
The neighbors of $z$ are similarly the neighbors for any
$\C_i$ containing $z$.

\section{Implementation}\label{s-implementation}
We have implemented the NCVX Python package for modeling
problems of the form~(\ref{e-prob}) and applying the NC-ADMM heuristic,
along with the relax-round-polish and relax methods.
The NCVX package is an extension of CVXPY \cite{cvxpy}.
The problem
objective and convex constraints are expressed using standard CVXPY
semantics. Nonconvex constraints are expressed implicitly by creating
a variable constrained to lie in one of the sets described in \S\ref{s-set}.
For example, the code snippet
\begin{verbatim}
x = Boolean()
\end{verbatim}
creates a variable $x \in \reals$ with the
implicit nonconvex constraint $x \in \{0,1\}$.
The convex relaxation, in this case $x \in [0,1]$, is
also implicit in the variable
definition.
The source code for NCVX is available at
\url{https://github.com/cvxgrp/ncvx}.

\subsection{Variable constructors}\label{create-var-sec}

The NCVX package provides the following functions for creating
variables with implicit nonconvex constraints,
along with many others not listed:
\begin{itemize}
\item \verb+Boolean(n)+ creates a variable
$x \in \reals^{n}$ with the implicit constraint
$x \in \{0,1\}^{n}$.
\item \verb+Integer(n, M)+ creates a variable
$x \in \reals^{n}$ with the implicit constraints
$x \in \integers^{n}$ and $\|x\|_\infty \leq \lfloor M \rfloor$.
\item \verb+Card(n, k, M)+ creates a variable
$x \in \reals^{n}$ with the implicit constraints
that at most $k$ entries are nonzero and $\|x\|_\infty \leq M$.
\item \verb+Choose(n, k)+ creates a variable
$x \in \reals^{n}$ with the implicit constraints that
$x \in \{0,1\}^{n}$ and has exactly $k$ nonzero entries.
\item \verb+Rank(m, n, k, M)+ creates a variable
$X \in \reals^{m \times n}$ with the implicit constraints
$\Rank(X) \leq k$ and $\|X\|_2 \leq M$.
\item \verb+Assign(m, n)+ creates a variable
$X \in \reals^{m \times n}$ with the implicit constraint
that $X$ is an assignment matrix.
\item \verb+Permute(n)+ creates a variable
$X \in \reals^{n \times n}$ with the implicit constraint
that $X$ is a permutation matrix.
\item \verb+Cycle(n)+ creates a variable
$X \in \reals^{n \times n}$ with the implicit constraint
that $X$ is the adjacency matrix of a Hamiltonian cycle.
\item \verb+Annulus(n,r,R)+ creates a variable $x \in \reals^n$ with
the implicit constraint that $r \leq \|x\|_2 \leq R$.
\item \verb+Sphere(n, r)+ creates a variable $x \in \reals^n$ with
the implicit constraint that $\|x\|_2 = r$.
\end{itemize}

\subsection{Variable methods}
Additionally, each variable created by the functions in \S\ref{create-var-sec}
supports the following methods:
\begin{itemize}
\item \verb+variable.relax()+ returns a list of convex constraints
that represent a convex relaxation of the nonconvex set $\C$,
to which the variable belongs.
\item \verb+variable.project(z)+ returns the Euclidean (or approximate)
projection of $z$ onto the nonconvex set $\C$,
to which the variable belongs.
\item \verb+variable.restrict(z)+ returns a list of convex constraints
describing the convex restriction $\Crestrict(z)$ at $z$ of the
nonconvex set $\C$, to which the variable belongs.
\item \verb+variable.neighbors(z)+ returns a list of neighbors $\Cneighbor(z)$
of $z$ contained in the
nonconvex set $\C$, to which the variable belongs.
\end{itemize}
Users can add support for additional nonconvex sets by implementing
functions that return variables with these four methods.

\subsection{Constructing and solving problems}
To construct a problem of the form (\ref{e-prob}), the user creates
variables $z_1,\ldots,z_k$ with the implicit constraints
$z_1 \in \C_1,\ldots,z_k \in \C_k$, where $\C_1,\ldots,\C_k$ are
nonconvex sets,
using the functions described in \S\ref{create-var-sec}.
The variable $z$ in problem (\ref{e-prob}) corresponds to the vector
$(z_1,\ldots,z_k)$.
The components of the variable $x$, the objective,
and the constraints are constructed using standard CVXPY syntax.

Once the user has constructed a problem object,
they can apply the following solve methods:
\begin{itemize}
\item \verb+problem.solve(method="relax")+ solves the convex relaxation
of the problem.
\item \verb+problem.solve(method="relax-round-polish")+
applies the relax-round-polish heuristic.
Additional arguments can be used to specify the parameters $K$ and $\lambda$.
By default the parameter values are $K=1$ and $\lambda = 10^4$.
When $K > 1$, the first sample $w_1 \in \reals^q$ is always 0.
Subsequent samples are drawn i.i.d.\ from
$N(0,\sigma^2I)$, where $\sigma$ is another parameter the user can set.
\item \verb+problem.solve(method="nc-admm")+
applies the NC-ADMM heuristic.
Additional arguments can be used to specify the number of starting points,
the number of iterations the algorithm is run from each starting point,
and the values of the parameters $\rho$, $\sigma$, and $\lambda$.
By default the algorithm is run from 5 starting points for 50 iterations,
the value of $\rho$ is drawn uniformly from $[0,1]$, and the other
parameter values are $\sigma = 1$ and $\lambda = 10^4$.
The first starting point is always
$z^0 = 0$ and subsequent starting points are drawn i.i.d.\ from
$\mathcal N(0,\sigma^2I)$.
\end{itemize}
The relax-round-polish and NC-ADMM methods
record the best point found $(x_\mathrm{best},z_\mathrm{best})$
according to the merit function.
The methods return the objective value $f_0(x_\mathrm{best},z_\mathrm{best})$ and
the residual $r(x_\mathrm{best},z_\mathrm{best})$,
and set the \verb+value+
field of each variable to the appropriate segment
of $x_\mathrm{best}$ and $z_\mathrm{best}$.

For example, consider the \emph{regressor selection} problem, which we will
discuss in \S\ref{regressor_prob}. This problem can be formulated as
\BEQ
\begin{array}{ll}
\mbox{minimize}   & \|Ax-b\|_2^2 \\
\mbox{subject to} & \|x\|_\infty \leq M \\
& \card(x)\leq k,
\end{array}
\label{e-card-ls}
\EEQ
with decision variable $x\in\reals^n$ and problem data $A\in\reals^{m\times n}$,
$b\in\reals^m$, $M > 0$, and $k \in \integers_+$. The following code attempts to approximately solve
this problem using our heuristic.

\begin{verbatim}
x = Card(n,k,M)
prob = Problem(Minimize(sum_squares(A*x-b)))
objective, residual = prob.solve(method="nc-admm")
\end{verbatim}
The first line constructs a variable
$x \in \reals^{n}$ with the implicit constraints
that at most $k$ entries are nonzero, $\|x\|_\infty \leq M$,
and $\|x\|_1 \leq kM$.
The second line
creates a minimization problem with objective $\|Ax-b\|_2^2$ and
no constraints.
The last line applies the NC-ADMM heuristic to the problem
and returns the objective value and residual of the best point found.

\section{Examples}\label{s-examples}
In this section we apply the NC-ADMM heuristic to a wide variety of hard
problems, \ie, that generally cannot be solved in polynomial time.
Extensive research has been done on specialized algorithms
for each of the problems discussed in this section.
Our intention is not to seek better performance than these specialized
algorithms, but rather to show that our
general purpose heuristic can yield decent results with minimal tuning.
Unless otherwise specified, the algorithm parameters are the defaults
described in \S\ref{s-implementation}. Whenever possible, we compare our heuristic to
GUROBI \cite{gurobi}, a commercial global optimization solver.
Since our implementation of
NC-ADMM supports minimal parallelization, we compare the number of
convex subproblems solved (and not the solve time).

\subsection{Regressor selection}\label{regressor_prob}
We consider the problem of approximating
a vector $b$ with a linear combination of at most $k$ columns of $A$
with bounded coefficients.
This problem can be formulated as
\BEQ
\begin{array}{ll}
\mbox{minimize}   & \|Ax-b\|_2^2 \\
\mbox{subject to} & \card(x)\leq k, \quad \|x\|_\infty \leq M,
\end{array}
\label{e-card-ls}
\EEQ
with decision variable $x\in\reals^n$ and problem data $A\in\reals^{m\times n}$,
$b\in\reals^m$, $k \in \integers_+$, and $M > 0$.
Lasso (least absolute shrinkage and selection operator) is
a well-known heuristic for solving this problem by adding $\ell_1$
regularization and minimizing $ \|Ax-b\|_2^2 +\lambda\|x\|_1$.
The value of $\lambda$ is chosen as the smallest value possible such
that $\card(x) \leq k$. (See
\cite[\S 3.4]{friedman2001elements} and
\cite[\S 6.3]{boyd2004convex}.)

\paragraph{Problem instances.}
We generated the matrix $A \in \reals^{m \times 2m}$
with i.i.d.\ $\mathcal N(0,1)$ entries, and
chose $b=A\hat x + v$, where $\hat x$ was drawn uniformly at random
from the set of vectors satisfying
$\card(\hat x) \leq \lfloor m/5 \rfloor$
and $\|x\|_\infty \leq 1$, and
$v \in \reals^m$ was a noise vector drawn from
$\mathcal{N}(0,\sigma^2 I)$.
We set $\sigma^2 = \|A \hat{x}\|^2/(400m)$ so that the signal-to-noise
ratio was near 20.

\begin{figure}
\begin{center}
\includegraphics[width=0.7\textwidth]{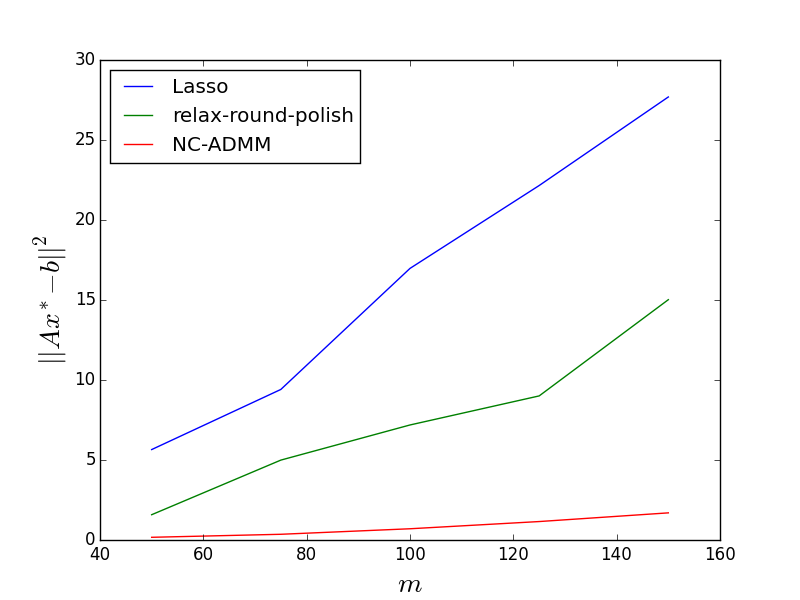}
\end{center}
\caption{The average error of solutions found by Lasso, relax-round-polish,
and NC-ADMM for 40 instances of the regressor selection problem.
}\label{regressor_results}
\end{figure}

\paragraph{Results.}
For each value of $m$, we generated $40$ instances of the problem
as described in the previous paragraph.
Figure \ref{regressor_results} compares the average sum of squares error
for the $x^*$ values found by the Lasso heuristic, relax-round-polish,
and NC-ADMM. For Lasso, we solved the problem for 100 values of $\lambda$
and then solved the polishing problem after fixing the sparsity pattern suggested by Lasso.

\subsection{3-satisfiability}
Given Boolean variables $x_1,\cdots, x_n$, a \emph{literal} is either a
variable or the negation of a variable, for example $x_1$ and $\neg x_2$.
A \emph{clause} is disjunction of literals (or a single literal), for
example $(\neg x_1 \lor x_2 \lor \neg x_3)$. Finally a formula is
in conjunctive normal form (CNF) if it is a conduction of clauses (or a single
clause), for example $(\neg x_1 \lor x_2 \lor \neg x_3) \land(x_1 \lor \neg x_2)$.
Determining the satisfiability of a formula in conjunctive normal form where each
clause is limited to at most three literals is called \emph{3-satisfiability} or
simply the \emph{3-SAT} problem. It is known
that 3-SAT is NP-complete, hence we do not expect to be able to solve a 3-SAT
in general using our heuristic.
A 3-SAT problem
can be formulated as the following
\BEQ
\label{3satlp}
\begin{array}{ll}
\mbox{minimize} &0\\
\mbox{subject to} & Az \leq b,\\
& z \in \{0,1\}^n,\\
\end{array}
\EEQ
where entries of $A$ are given by
\[
a_{ij} =
\begin{cases}
-1 & \mbox{if clause $i$ contains $x_j$}\\
1 & \mbox{if clause $i$ contains $\neg x_j$}\\
0 & \mbox{otherwise},
\end{cases}
\]
and the entries of $b$ are given by
\[
b_{i} = (\mbox{number of negated literals in clause $i$}) -1.
\]

\paragraph{Problem instances.}
We generated 3-SAT problems with varying numbers of clauses and variables
randomly as in \cite{mitchell1992hard,lipp2014variations}.
As discussed in \cite{crawford1996experimental},
there is a threshold around $4.25$ clauses
per variable when problems transition from being feasible to being infeasible.
Problems near this threshold are generally found
to be hard satisfiability problems.
We generated 10 instances for each choice of number of clauses and variables,
verifying that each instance is feasible using GUROBI \cite{gurobi}.

\begin{figure}
\begin{center}
\includegraphics[width=0.7\textwidth]{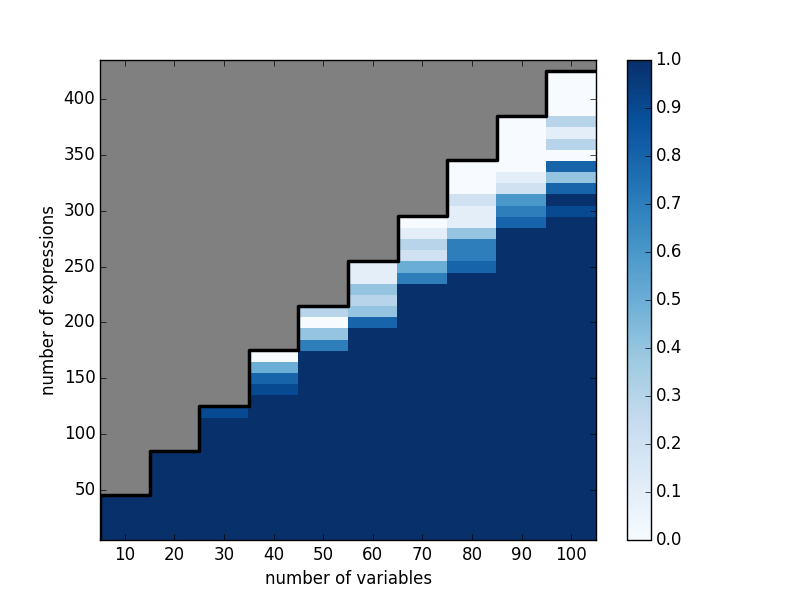}
\end{center}
\caption{Fraction of runs for which a satisfying assignment to
random 3-SAT problems were found for problems of varying sizes.
The problems in the gray region were not tested.}
\label{3sat_results}
\end{figure}

\paragraph{Results.}
We ran NC-ADMM heuristic on each instance, with 10 restarts,
and 100 iterations, and we chose the step size $\rho=10$.
Figure~\ref{3sat_results} shows the fraction of instances solved correctly
with NC-ADMM. We see that using this heuristic, satisfying assignments
can be found consistently for up to 3.2
constraints per variable, at which point success starts to decrease.
Problems in the gray region in figure~\ref{3sat_results} were not
tested since they are infeasible with high probability.
We also tried the relax-round-polish heuristic, but it often failed to
solve problems with more than $50$ clauses.

\begin{incomplete}
\subsection{Graph coloring}
Given a graph $G$, the smallest number of colors needed to color the vertices of
$G$ such that no two adjacent vertices have the same color is called the
\emph{chromatic} number of $G$ and is often denoted by $\chi(G)$. It is easy to
show that $\omega(G)\leq\chi(G)\leq \Delta(G) + 1$,
where $\omega(G)$ is the size of a maximum clique in $G$ and
$\Delta(G)$ is the maximum degree in $G$.
Finding the chromatic number of a graph is known to be NP-hard.
For a given graph $G$ we can construct a \emph{coloring matrix}
$Z\in\{0,1\}^{n\times \Delta(G) + 1}$,
where $Z_{ij}=1$ means that vertex $i$ is colored with color $j$.
Each vertex must be colored with exactly one color, so there should be
exactly one nonzero entry in every row. In other words, $Z\ones = \ones$.
For any adjacent vertices $i$ and $j$,
we must have $Z_{i,:}+Z_{j,:}\leq\ones$,
where $Z_{i,:}$ denotes the $i$th row of $Z$.
We want to minimize the number of colors used,
which can be written as $\sum_j \| Z_{:,j}\|_\infty$, since
$\| Z_{:,j}\|_\infty$ is $0$ if color $j$ is not used, and $1$ if it is.
If $Z^\star$ is an optimal solution, permuting the columns of $Z^\star$
gives additional optimal solutions.
To break this symmetry, we add the penalty term
\[
\tau(Z) = \sum_{i=1}^{n-1} \max\{h^T(Z_{:,i} - Z_{:,i+1}),0\},
\]
where $h \in \reals^n$ is an arbitrary vector.
Adding the penalty term does not alter the problem, since there
is always a permutation of $Z^\star$ for which $\tau(Z) = 0$.
The graph coloring problem can thus be formulated as
\BEQ
\label{graphcolor}
\begin{array}{ll}
\mbox{minimize} &\| Z_{:,j}\|_\infty + \tau(Z)\\
\mbox{subject to} & Z_{i,:}+Z_{j,:}\leq \ones \quad \mbox{$i,j$ adjacent}\\
& Z_{i,:} \in \mathcal{C}, \quad i=1,\ldots,n,
\end{array}
\EEQ
where $\C = \{ x\in \{0,1\}^n \mid \card(x) = 1\}$,
or \verb+Choose(n,1)+ in CVXPY.

\paragraph{Problem instances.}
We generated random graphs drawn uniformly from the set of graphs
with $n$ vertices and $\lfloor n(n-1)/2k \rfloor$ edges,
for a range of $n$ and $k$.
The parameter $k$ is the ratio of edges in a complete graph to
edges in the generated graph.
We generated 10 graphs for each choice of $n$ and $k$.
We found the chromatic number of each graph by solving problem
(\ref{graphcolor}) using GUROBI.
The vector $h \in \reals^n$ for each problem instance was chosen
uniformly at random on the $n$-sphere.

\paragraph{Results.}
Figure \ref{graph_coloring_results} shows the average difference between
the result found by NC-ADMM and the true chromatic number.
NC-ADMM was run from 5 random starting points for 50 iterations.
The value of $\rho$ for each starting point was drawn uniformly at
random from $[0,1]$.

\begin{figure}
\begin{center}
\includegraphics[width=0.7\textwidth]{graph_coloring_results}
\end{center}
\caption{Graph coloring results.
}\label{graph_coloring_results}
\end{figure}
\end{incomplete}

\begin{incomplete}
\subsection{Assignment (not final yet)}
In the most general form of \emph{assignment problem}, there are $n$ tasks and
$m$ processors. Processor $j$ can perform task $i$ with cost $W_{ij}$. The problem
is assigning tasks to processors so that the overall cost is minimized assuming that
each processor can be assigned to at most one task. This problem is equivalent to
finding a minimum weight matching of a bipartite graph

The goal is to find an assignment matrix $X\in\{0,1\}^{m\times n}$ such that
$\Tr(X^TW)$ is minimized. So the problem can be formulated as
\BEQ
\label{assignment}
\begin{array}{ll}
\mbox{minimize} & \Tr(Z^TW)\\
\mbox{subject to} & Z \quad \mbox{assignment matrix},\\
\end{array}
\EEQ
\end{incomplete}

\subsection{Circle packing}
In \emph{circle packing problem} we are interested in finding the smallest square
in which we can place $n$ non-overlapping circles with radii $r_1,\ldots,r_n$
\cite{goldberg1970packing}.
This problem has been studied extensively
\cite{stephenson2005introduction, castillo2008solving, collins2003circle}
and a database of densest known packings for different numbers of circles can be found
in \cite{Specht2013}.
The problem can be formulated as
\BEQ
\label{packing}
\begin{array}{ll}
\mbox{minimize} & l\\
\mbox{subject to} & r_i\ones \leq x_i \leq (l-r_i) \ones, \quad i = 1,\ldots,n\\
& x_i-x_j = z_{ij},\quad i= 1,\ldots , n-1,\quad j = i+1,\ldots, n \\
& 2\sum_{k=1}^n r_i \geq \|z_{ij}\|_2 \geq r_i + r_j, \quad i= 1,\ldots , n -1,\quad j = i+1,\ldots, n,
\end{array}
\EEQ
where $x_1,\ldots,x_n \in \reals^2$ are variables representing the circle centers
and $z_{12},z_{13},\ldots,z_{n-1,n} \in \reals^2$ are additional variables
representing the offset between pairs $(x_i,x_j)$.
Note that each $z_{ij}$ is an element of an annulus.

\paragraph{Problem instances.}
We generated problems with different numbers of circles. Here we report the performance
of the relax-round-polish heuristic for a problem with $n=41$,
in two cases:
a problem with all circle radii equal to $0.5$, and a problem where the radii
were chosen uniformly at random from the interval $[0.2,0.5]$.

\begin{figure}
\begin{center}
\includegraphics[width=.46\linewidth]{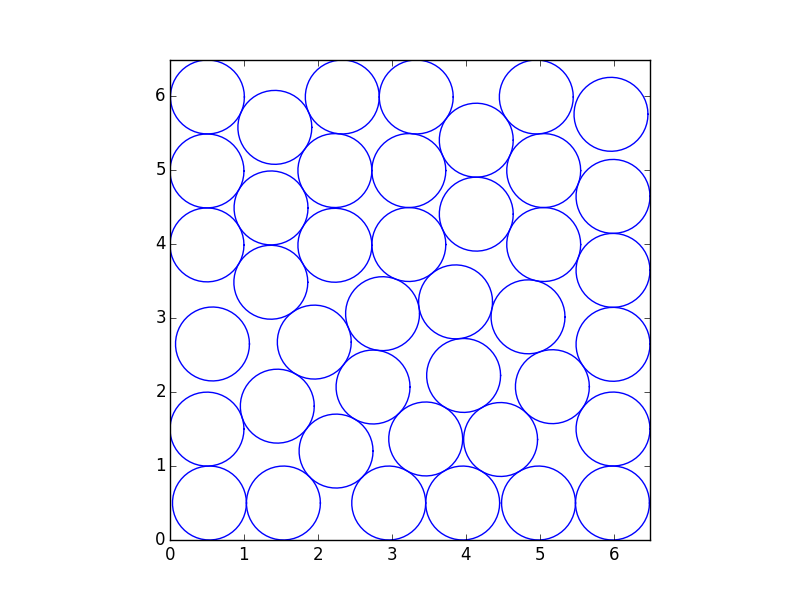}
\hspace*{\fill}
\includegraphics[width=.46\linewidth]{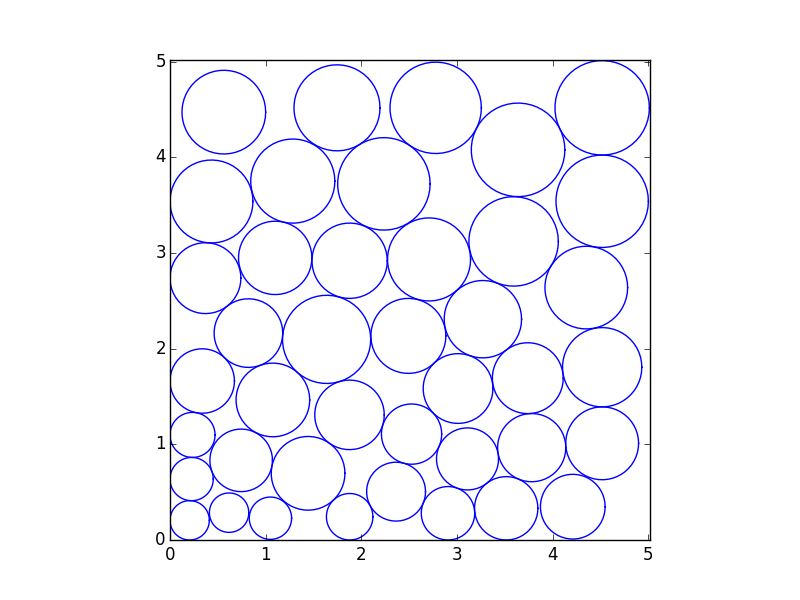}
\caption{
Packing for $n=41$ circles with equal and different radii.
\label{circles41}}
\end{center}
\end{figure}

\paragraph{Results.}
We run the relax-round-polish heuristic in both cases.
For this problem, the heuristic is effectively equivalent to
many well-known methods like the convex-concave procedure and the majorization-minimization
(MM) algorithm. Figure \ref{circles41} shows the packing found by our heuristic for $n=41$.
The obtained packing covers $78.68\%$ of the area of the bounding square, which is close
to the densest known packing, which covers $79.27\%$ of the area.
We observed that NC-ADMM is no more effective than
relax-round-polish for this problem.

\subsection{Traveling salesman problem}
In the traveling salesman problem (TSP), we wish to find
the minimum weight Hamiltonian cycle in a weighted graph.
A Hamiltonian cycle is a path that starts and ends on the same
vertex and visits each other vertex in the graph exactly once.
Let $G$ be a graph with $n$ vertices and $D\in\symm^n$ be the (weighted)
adjacency matrix, \ie, the real number $d_{ij}$ denotes the distance between $i$ and $j$.
We can formulate the TSP problem for $G$ as follows
\BEQ
\label{tsp}
\begin{array}{ll}
\mbox{minimize} & (1/2)\Tr(D^TZ)\\
\mbox{subject to} &Z \in \mathcal H_n,
\end{array}
\EEQ
where $Z$ is the decision variable \cite{lawler1985traveling,kruskal1956shortest,
dantzig1954solution,hoffman2013traveling}.

\paragraph{Problem instances.}
We generated $n=75$ points in $[-1,1]^2$.
We set $d_{ij}$ to be the Euclidean distance between
points $i$ and $j$.

\paragraph{Results.}

Figure \ref{tsp} compares the Hamiltonian cycle
found by the NC-ADMM heuristic, which had cost $14.47$,
with the optimal Hamiltonian cycle,
which had cost $14.16$.
The cycle found by our heuristic has a few clearly suboptimal paths,
but overall is a reasonable approximate solution.
We ran NC-ADMM with 5 restarts and 100 iterations.
GUROBI solved $4190$ subproblems before finding a solution as
good as that found by NC-ADMM,
which solved only 500 subproblems.
The relax-round-polish heuristic does not perform well on this problem.
The best objective value found by the heuristic is $35.6$.

\begin{figure}
\begin{center}
\begin{psfrags}
\psfrag{x}[B][B]{\raisebox{-1.2ex}{\tiny $t$}}
\psfrag{p1}[B][B]{\raisebox{+.5ex}{\tiny{$P^\mathrm{Eng}_t$}}}
\psfrag{p2}[B][B]{\raisebox{+.5ex}{{\tiny $P^\mathrm{batt}_t$}}}
\psfrag{p3}[B][B]{\raisebox{+.5ex}{\tiny{$E_t$}}}
\psfrag{p4}[B][B]{\raisebox{+.5ex}{\tiny{$z_t$}}}
\includegraphics[width=.46\linewidth]{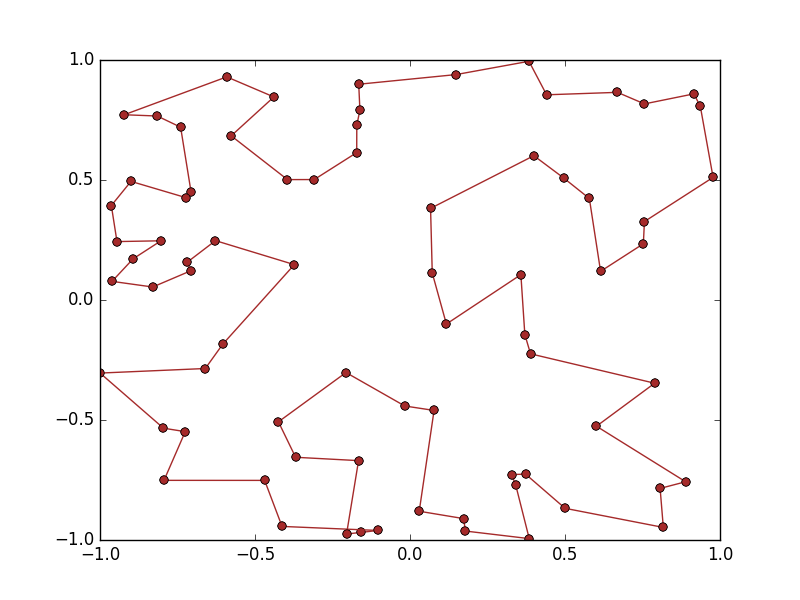}
\end{psfrags}
\hspace*{\fill}
\begin{psfrags}
\psfrag{x}[B][B]{\raisebox{-1.2ex}{\tiny $t$}}
\psfrag{p1}[B][B]{\raisebox{+.5ex}{\tiny{$P^\mathrm{Eng}_t$}}}
\psfrag{p2}[B][B]{\raisebox{+.5ex}{{\tiny $P^\mathrm{batt}_t$}}}
\psfrag{p3}[B][B]{\raisebox{+.5ex}{\tiny{$E_t$}}}
\psfrag{p4}[B][B]{\raisebox{+.5ex}{\tiny{$z_t$}}}
\includegraphics[width=.46\linewidth]{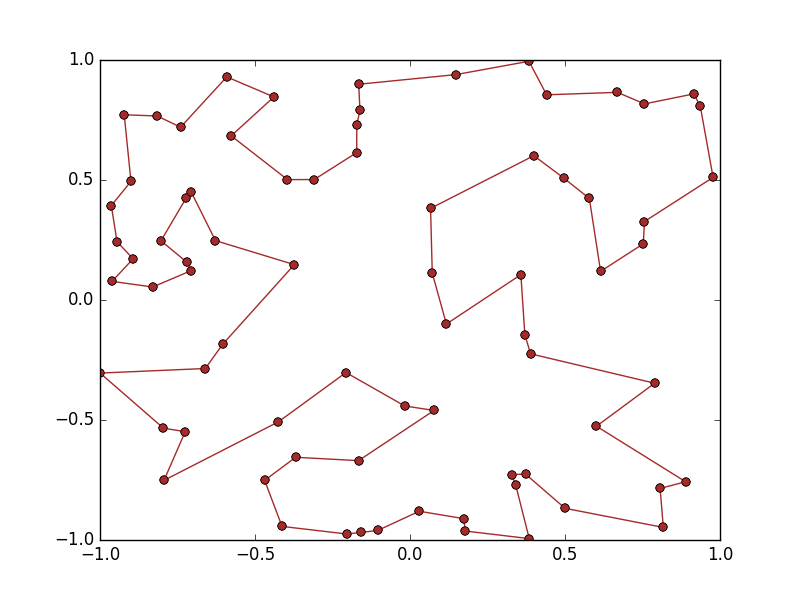}
\end{psfrags}
\caption{
Left: Hamiltonian cycle found by NC-ADMM.
Right: optimal Hamiltonian cycle, found using GUROBI.
\label{tsp}}
\end{center}
\end{figure}

\subsection{Factor analysis model}
The factor analysis problem decomposes a matrix as a sum of a low-rank
and a diagonal matrix and has been studied extensively (for example
in \cite{saunderson2012diagonal,ning2015linear}).
It is also known as the \emph{Frisch} scheme in
the system identification literature \cite{kalman1985identification,
david1993opposite}.
The problem is the following
\BEQ
\label{factor}
\begin{array}{ll}
\mbox{minimize} & \|\Sigma-\Sigma^\mathrm{lr}-D\|_F^2\\
\mbox{subject to} & D = \diag(d),\quad d\geq 0\\
&\Sigma^\mathrm{lr} \succeq 0\\
&\Rank(\Sigma^\mathrm{lr}) \leq k,

\end{array}
\EEQ
where $\Sigma^\mathrm{lr}\in\symm_+^{n}$ and diagonal matrix
$D\in\reals^{n\times n}$ with nonnegative diagonal entries
are the decision variables, and $\Sigma\in\symm_+^n$ and
$k\in\integers_+$ are problem data.
One well-known heuristic for solving this problem is adding
$\|\cdot\|_*$, or nuclear norm, regularization and minimizing
$ \|\Sigma-\Sigma^\mathrm{lr}-D\|_F^2 +\lambda\|\Sigma^\mathrm{lr}\|_*$.
The value of $\lambda$ is chosen as the smallest value possible such
that $\Rank(\Sigma^\mathrm{lr}) \leq k$.
Since $\Sigma^\mathrm{lr}$ is positive semidefinite,
$\|\Sigma^\mathrm{lr}\|_* = \Tr(\Sigma^\mathrm{lr})$.

\paragraph{Problem instances.}
We set $k= \lfloor n/2 \rfloor$ and
generated the matrix $F \in \reals^{n \times k}$
by drawing the entries i.i.d.\ from a standard normal distribution.
We generated a diagonal matrix $\hat D$ with
diagonal entries drawn i.i.d.\ from an exponential distribution
with mean $1$.
We set $\Sigma = FF^T + \hat D + V$,
where $V \in \reals^{n \times n}$ is a noise matrix with
entries drawn i.i.d.\ from $\mathcal{N}(0,\sigma^2)$.
We set $\sigma^2 = \|FF^T + \diag(d)\|^2_F/(400n^2)$ so that
the signal-to-noise ratio was near 20.

\begin{figure}
\begin{center}
\includegraphics[width=0.7\textwidth]{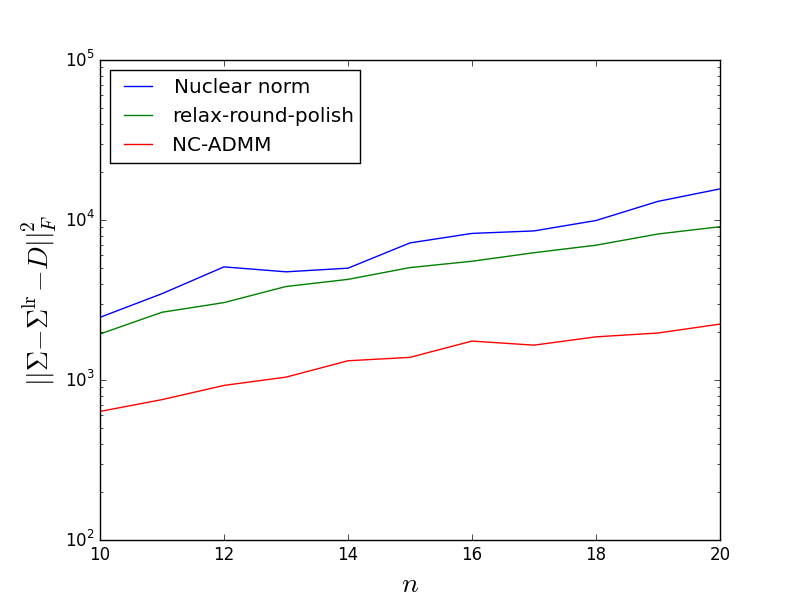}
\end{center}
\caption{The average error of solutions found by the nuclear norm,
relax-round-polish, and NC-ADMM heuristics for 50 instances of the factor analysis problem.
}\label{factor_analysis_results}
\end{figure}

\paragraph{Results.}
Figure \ref{factor_analysis_results} compares the average sum of squares error
for the $\Sigma^\mathrm{lr}$ and $d$ values found by NC-ADMM,
relax-round-polish, and the nuclear norm heuristic over 50 instances per value of $n$.
We observe that the sum of squares error obtained by NC-ADMM is smaller
than that obtained by the nuclear norm and relax-round-polish heuristics.

\subsection{Job selection}
In the job selection problem there are $n$ jobs and $m$ resources.
Each job $i$ consumes $A_{ij} \geq 0$ units of resource $j$, and
up to $d_i > 0$ instances of job $i$ can be accepted.
Executing job $i$ produces $c_i > 0$ units of profit.
The goal is to maximize profit subject to the constraint that
at most $b_j > 0$ units of each resource $j$ are consumed.
The job selection problem can be formulated as
\BEQ
\label{job}
\begin{array}{ll}
\mbox{maximize} & c^Tz \\x
\mbox{subject to} & Az \leq b \\
& 0 \leq z \leq d\\
& z \in \integers^n,
\end{array}
\EEQ
where $z$ is the decision variable and $A \in \reals^{m \times n}$,
$b \in \reals^m_+$, and $d \in \integers_+^n$ are problem data.
This problem is NP-hard in general.
When $m=1$, this problem is
equivalent to the \emph{knapsack problem}, which
has been studied extensively; see, \eg,
\cite{chu1998genetic, chekuri2005polynomial}.

\paragraph{Problem instances.}
We set $m=\lfloor n/10 \rfloor$ and generated $A \in \reals^{m \times n}$
by randomly selecting $\lfloor mn/10 \rfloor$ entries to be nonzero.
The nonzero entries were drawn i.i.d.\ from the uniform distribution
over $[0,5]$.
Entries of $c \in \reals^n$ were drawn i.i.d.\ from
the uniform distribution over $[0,1]$.
Entries of $d \in \integers^n$ were drawn i.i.d.\ from the uniform
distribution over the set $\{1,\ldots,5\}$.
We generated $b \in \reals^{m}$ by first generating
$\hat{z} \in \integers^n$, where each $\hat{z}_i$ was drawn
from the uniform distribution over the set $\{0,\ldots,d_i\}$,
and then setting $b = A\hat{z}$.

\begin{figure}
\begin{center}
\includegraphics[width=0.7\textwidth]{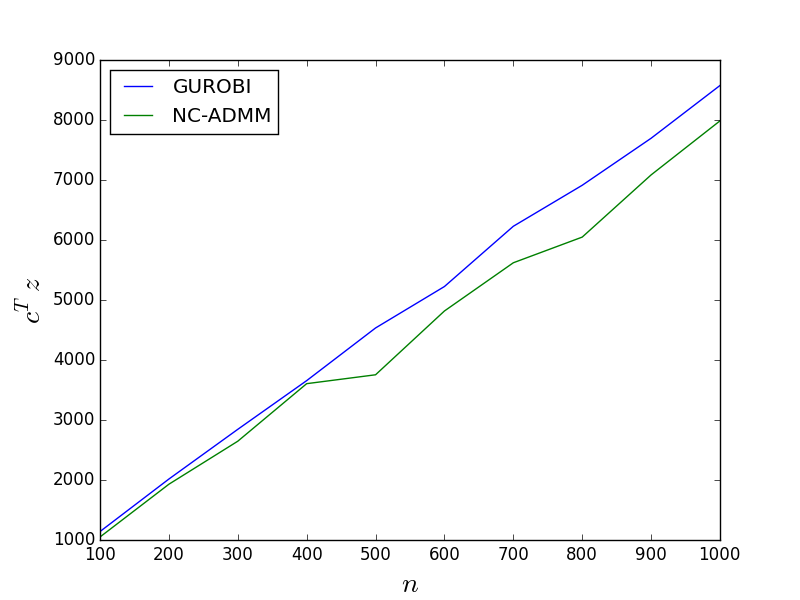}
\end{center}
\caption{The average objective value for the job selection problem over 10 different instances.
}\label{job_selection_results}
\end{figure}

\paragraph{Results.}
We generated problem instances for a range of $10 \leq n \leq 100$.
Figure \ref{job_selection_results} compares, for each $n$,
the average value of $c^Tz$ found by the NC-ADMM heuristic and by GUROBI over $10$ instances.
NC-ADMM was run from 10 random starting points for 100 iterations.
The value of $\rho$ for each starting point was drawn
from the uniform distribution over $[0,5]$.
GUROBI's run time was limited to 10 minutes.
NC-ADMM always found a feasible $z$ with an objective value not much
worse than that found by GUROBI.
We also tried the relax-round-polish heuristic on the problem instances,
but it never found a feasible $z$.

\begin{incomplete}
\subsection{Interpretable models}
Developing interpretable patient-level predictive models using medical data has
attracted researchers attention \cite{LethamRuMcMa15}. A widely used example
of an interpretable model is CHADS\textsubscript2 \cite{gage2001validation}
to predict stroke in patients with atrial fibrillation,
in which the patient's score is calculated by assigning one
or two `points'
for the presence of $5$ different features (such as age above $75$).

In a similar framework here we study a logistic regression model
to predict a binary result (\eg, stroke or no stroke),
where the
coefficients (individual scores) come from a discrete
set for example $\{0,\pm1,\pm2\}$. Let
$u_1,\ldots, u_q\in \reals^n$ denote the training points where the
result is positive and let $u_{q+1},\ldots, u_m\in \reals^n$ denote the
training points where the result is negative. The problem is
\BEQ
\begin{array}{ll}
\mbox{minimize}   & -\sum_{i=1}^q(a^Tu_i + b) + \sum_{i=1}^n
\log\left(1+\exp(a^Tu_i+b)\right)\\
\mbox{subject to} & a \in \{0,\pm1,\pm2\}^n,
\end{array}
\label{e-card-ls}
\EEQ
with decision variables $a,b$.

\paragraph{Problem instances.}
TODO

\paragraph{Results.}
TODO
\end{incomplete}

\begin{incomplete}
\subsection{Pooling problems}
The pooling problem is a multi-commodity flow problem, where each node
is representative of a pool where (typically fluid) material are continuously mixed
\cite{gupte2013pooling, bodington1990history}.
Consider a directed acyclic graph with $n$ nodes where each edge $e$
is characterized by a vector $m_e\in\reals_+^q$
representing the flow of each of
$q$ constituents.
For each node $v$, let $D_v^+$($D_v^-$) be the set of edges entering to
(exiting from) node
$v$. At each node
$v$, we have the (mass conservation) equations
\[
\sum_{e\in D_v^+}m_e = \sum_{e\in D_v^-}m_e.
\]
Also, the mass flows from any node should have the same concentration of each of
$q$ constituents. In other words, for any node $v$ there exists a vector $c_v$ such
that for every edge $e\in D_v^-$,
\[
c_v = \frac{m_e}{\ones^Tm_e}.
\]
Some nodes are \emph{source} nodes, with a single input edge with known
mass flow rates. Some vertices are \emph{sink} nodes, with a single output
edge with desired mass flow rates.
Let $I_+$($I_i$) denote the set of source (sink) nodes.
In the simplest case, we pay $p_v$ units
for every unit of the raw material in node $v\in I_+$,
and get paid $p_w$ units for every unit of the blended material in node $w\in I_-$.
For nodes $v\in I_+$, the vector $c_v$ is given and there is a upper bound on
the total amount of input $\ones^T m_e \leq B_e$ for $e \in D_v^+$.
Hence the problem is
\BEQ
\label{pooling}
\begin{array}{ll}
\mbox{minimize} & \sum_{v\in I_+}p_v\sum_{e\in D_v^+}m_e -
\sum_{v\in I_-}p_v\sum_{e\in D_v^-}m_e\\
\mbox{subject to} & \sum_{e\in D_v^+}m_e = \sum_{e\in D_v^-}m_e
\quad v=1,\ldots,n\\
& c_v = \frac{m_e}{\ones^Tm_e} \quad v=1,\ldots,n, \quad e\in D_v^-\\
& \ones^T m_e \leq B_e \quad v \in I_-,\quad e\in D_v^+\\
& \tilde c_v \leq c_v \leq \tilde c_v \quad v \in I_+.
\end{array}
\EEQ
Notice that if a network only consists of source and sink nodes, this problem
is convex. Also for fixed concentration vectors $c_v$, this is a convex problem.

\paragraph{Problem instances.}
TODO the class isn't in the Python package yet.

\paragraph{Results.}
TODO
\end{incomplete}

\subsection{Maximum coverage problem}
A collection of sets $\mathcal S = \{S_1,S_2,\ldots, S_m\}$
is defined over a domain of elements $\{e_1,e_2,\ldots,e_n\}$ with associated weights
$\{w_i\}_{i=1}^n$. The goal is to find the collection of no more than $k$
sets $\mathcal S^\prime \subseteq \mathcal S$ that maximizes the total weight
of elements covered by $\mathcal S^\prime$
\cite{khuller1999budgeted,hochbaum1996approximating,cohen2008generalized}.
Let $x_i\in\{0,1\}$,  for  $i=1,\ldots, n$, be
a variable that takes $1$ if element $e_i$ is covered and $0$ otherwise.
Let $y \in\{0,1\}^m$ be
a variable with entry $y_j = 1$ if set $j$ is selected.
The problem is
\BEQ
\label{coverage}
\begin{array}{ll}
\mbox{maximize} & w^Tx \\
\mbox{subject to} & \sum_{j \in S_j}y_j \geq x_i, \quad i=1,\ldots,n\\
& x_i \in\{0,1\}, \quad i=1,\ldots,n\\
& y \in \{0,1\}^m \\
& \card(y) = k.
\end{array}
\EEQ
Note that $y$ is a Boolean vector with fixed cardinality.

\paragraph{Problem instances.}
We generated problems as follows.
Each set contained each of the elements independently
with a constant probability $p$. Hence the expected size of each set was $np$.
There were $m=3/p$ sets, so the expected total
number of elements in all sets (with repetition) was equal to $mnp=3n$.
We set $k=1/(3p)$.
Each $w_i$ was chosen uniformly at random from the interval $[0,1]$.

\begin{figure}
\begin{center}
\includegraphics[width=0.7\textwidth]{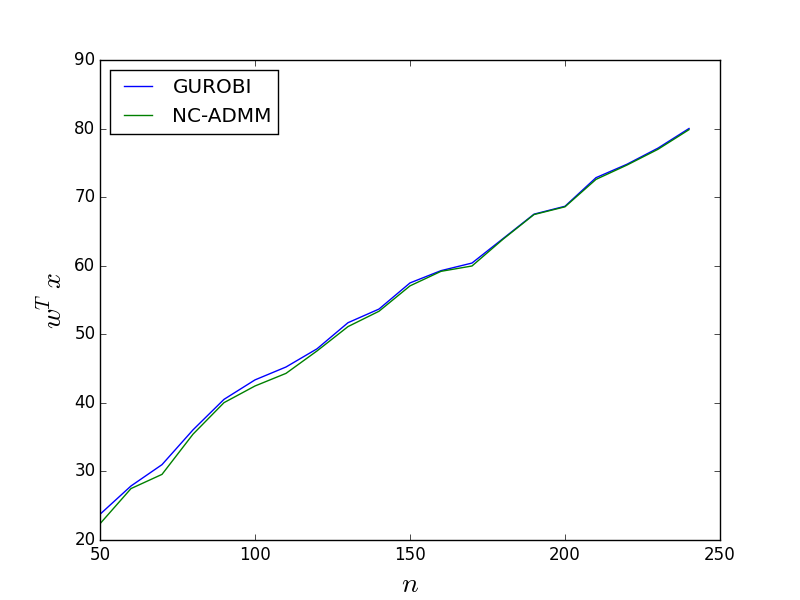}
\end{center}
\caption{The average solution weight over 10 different instances.
}\label{max_coverage_results}
\end{figure}

\paragraph{Results.}
We generated problems as described above for $n=50,60,\ldots,240$ and $p=0.01$.
For each value of $n$, we generated $10$ problems and recorded
the average weight $w^Tx$ of the approximate solutions found by NC-ADMM
and the optimal solutions found by GUROBI.
Figure \ref{max_coverage_results} shows the results of our comparison
of NC-ADMM and GUROBI.
Approximate solutions found by the relax-round-polish
heuristic were far worse than those found by NC-ADMM for this problem.

\subsection{Inexact graph isomorphism}
Two (undirected) graphs are isomorphic if we can permute the vertices of one so it
is the same as the other (i.e., the same pairs of vertices are connected by edges).
If we describe them by their adjacency matrices $A$ and $B$,
isomorphism is equivalent to
the existence of a permutation matrix $Z\in\reals^{n\times n}$
such that $ZAZ^T = B$, or equivalently $ZA=BZ$.

Since in practical applications isomorphic graphs might be contaminated
by noise, the inexact graph isomorphism problem is usually stated
\cite{aflalo2014graph, umeyama1988eigendecomposition,cross1997inexact},
in which we want to find a permutation matrix $Z$ such that the
disagreement $\|ZAZ^T-B\|_F^2$
between the transformed matrix and the target matrix is minimized.
Since $\|ZAZ^T-B\|_F^2=\|ZA-BZ\|_F^2$ for any permutation matrix $Z$,
the inexact graph isomorphism problem can be formulated as
\BEQ
\label{isomorphism}
\begin{array}{ll}
\mbox{minimize} & \|ZA-BZ\|_F^2\\
\mbox{subject to} & Z\in \mathcal P_n.
\end{array}
\EEQ
If the optimal value of this problem is zero, it means that $A$ and $B$
are isomorphic. Otherwise, the solution of this problem minimizes the
disagreement of $ZAZ^T$ and $B$ in the Frobenius norm sense.

Solving inexact graph isomorphism problems is of interest in pattern recognition
\cite{conte2004thirty, rocha1994shape}, computer vision
\cite{schellewald2001evaluation}, shape analysis
\cite{sebastian2004recognition, he2006object},
image and video indexing \cite{lee2006graph}, and neuroscience
\cite{vogelstein2011large}.
In many of the aforementioned fields
graphs are used to represent geometric structures, and
$\|ZAZ^T-B\|_F^2$
can be interpreted as the strength of geometric deformation.

\paragraph{Problem instances.}
It can be shown that if $A$ and $B$ are isomorphic and
$A$ has distinct eigenvalues
and for all
eigenvectors $v$ of $A$ for which $\ones^T v \neq 0$, then the relaxed
problem has a unique solution which is the permutation matrix that
relates $A$ and $B$ \cite{aflalo2014graph}.
Hence, in order to generate harder problems,
we generated the matrix $A$ such that it violated these conditions.
In particular, we constructed $A$ for the Peterson graph ($3$-regular with $10$ vertices),
icosahedral graph ($5$-regular with $12$ vertices),
Ramsey graph ($8$-regular with $17$ vertices),
dodecahedral graph ($3$-regular with $20$ vertices),
and the Tutte-Coxeter graph ($3$-regular with $30$ vertices).
For each example we randomly permuted the vertices to obtain two isomorphic graphs.

\paragraph{Results.}
We ran NC-ADMM with $20$ iterations and $5$ restarts.
For all of our examples NC-ADMM was able to find the permutation relating
the two graphs.
It is interesting to notice that running the algorithm multiple times
can find different solutions if there is more than one permutation
relating the two graphs.
The relax-round-polish heuristic failed to find a solution for all
of the aforementioned problems.

\begin{incomplete}
\subsection{Phase retrieval}
The problem of recovery of a signal given the magnitude of its Fourier
transform
\cite{candes2015phase,shechtman2014phase}
arises in various fields of science and engineering, including
X-ray crystallography \cite{harrison1993phase,millane1990phase},
microscopy \cite{miao2008extending},
astronomy \cite{fienup1987phase},
diffraction and array imaging \cite{bunk2007diffractive},
and optics \cite{walther1963question}.
In its most general form, the phase retrieval problem can be formulated as
finding $x\in\reals^n$ such that $|(Fx)_i|$ is given for $i=1,\ldots,m$,
where $F\in\complex^{m\times n}$. (A special case is when $Fx$ represents
the Fourier transform of $x$.) Decomposing $F$ into its real and imaginary
part by $F=F_R+iF_I$ for $F_R,F_I\in\reals^{m\times n}$, we have the following
formulation
\BEQ
\label{phase}
\begin{array}{ll}
\mbox{find} & x\\
\mbox{subject to} & F_Rx = u\\
& F_Ix=v\\
&u_i^2+v_i^2=y_i^2\quad i=1,\ldots,m,
\end{array}
\EEQ
with decision variable $x,u,v$ and problem data $F,y$.
\end{incomplete}

\begin{incomplete}
\subsection{Sudoku}
A Sudoku puzzle is a partially completed $N\times N$ grid of cells divided
into N segments. The goal is to fill the grid using a prescribed set of
$N$ distinct symbols such that each element of the set appears exactly
once in each row, column, and segment. We define $N^3$ binary variables
$z_{ijk}$ for $i,j,k=1,\ldots, n$ where $z_{ijk}=1$ if and only the element $(i,j)$
is filled with $k$. The problem then can be formulated as
\[
\label{sudoku}
\begin{array}{ll}
\mbox{find} & z_{ijk}\quad1\leq i,j,k \leq n\\
\mbox{subject to} &\sum_{k} z_{ijk} = 1 \quad 1\leq i,j \leq n\\
&\sum_{i} z_{ijk} = 1 \quad 1\leq j,k \leq n\\
&\sum_{j} z_{ijk} = 1 \quad 1\leq i,k \leq n\\
&\sum_{k} z_{ijk} = 1 \quad 1\leq i,j \leq n\\
&\sum_{(i,j)\in\mathcal I} z_{ijk} = 1 \quad \mathcal I~\mbox{is a region}\\
&z_{ijk}=1 \quad A_{ij}=k
\end{array}
\]

It is shown that solving Sudoku puzzles is NP-complete
\cite{takayuki2003complexity,kendall2008survey}.
In \cite{Derbinsky2013AnImp} a message-passing version of ADMM has been used to solve Sudoku puzzles.
\end{incomplete}

\section{Conclusions}
We have discussed the relax-round-polish and NC-ADMM heuristics and
demonstrated their performance on many different problems with
convex objectives and decision variables from a nonconvex set.
Our heuristics are easy to extend to additional problems because
they rely on a simple mathematical interface for nonconvex
sets.
We need only know a method for (approximate) projection onto the set.
We do not require but benefit from knowing
a convex relaxation of the set, a convex restriction at any point in the set,
and the neighbors of any point in the set under some discrete distance metric.
Adapting our heuristics to any particular problem is straightforward,
and we have fully automated the process in the NCVX package.

We do not claim that our heuristics give state-of-the-art results
for any particular problem.
Rather, the purpose of our heuristics is to give a fast and reasonable
solution with minimal tuning for a wide variety of problems.
Our heuristics also take advantage of the tremendous progress in
technology for solving general convex optimization problems,
which makes it practical to treat solving a convex problem
as a black box.

\clearpage
\nocite{*}
\bibliography{noncvx_admm}

\end{document}